\documentclass[a4paper,11pt]{article}

\usepackage[T1]{fontenc}
\usepackage{lmodern}

\usepackage{dsfont}

\usepackage[utf8]{inputenc}
\usepackage{amsmath}
\usepackage{amsthm}
\usepackage{amssymb}
\usepackage{mathrsfs}
\usepackage{graphicx}
\usepackage[all]{xy}
\usepackage{hyperref}
\usepackage[usenames, dvipsnames]{xcolor}
\usepackage{soul}
\usepackage{thm-restate}

\usepackage{stmaryrd}
\usepackage{caption}
\usepackage{abstract} 
\usepackage[shortlabels]{enumitem}

\usepackage{tikz} 

\newtheorem{thm}{Theorem}[section]
\newtheorem{cor}[thm]{Corollary}
\newtheorem{claim}[thm]{Claim}

\newtheorem{lemma}[thm]{Lemma}
\newtheorem{prop}[thm]{Proposition}

\theoremstyle{definition}
\newtheorem{definition}[thm]{Definition}
\newtheorem{ex}[thm]{Example}
\newtheorem{remark}[thm]{Remark}

\newtheorem*{definition*}{Definition}

\newcommand{\hooklongrightarrow}{\lhook\joinrel\longrightarrow}

\DeclareMathOperator{\stab}{\mathrm{stab}}
\DeclareMathOperator{\Mod}{\mathrm{Mod}}
\DeclareMathOperator{\amod}{\mathfrak{mod}}

\DeclareMathOperator{\Z}{\mathbb{Z}}
\DeclareMathOperator{\id}{{id}}

\def\rquotient#1#2{%
	\makeatletter
	\raise.3ex\hbox{$#1$}/\lower.3ex\hbox{$#2$}%
	\makeatother
}	

\makeatletter
\newcommand{\subjclass}[2][2010]{%
	\let\@oldtitle\@title%
	\gdef\@title{\@oldtitle\footnotetext{#1 \emph{Mathematics subject classification.} #2}}%
}
\newcommand{\keywords}[1]{%
	\let\@@oldtitle\@title%
	\gdef\@title{\@@oldtitle\footnotetext{\emph{Key words and phrases.} #1.}}%
}
\makeatother

\newcommand{\Address}{{
		\bigskip
		\small
		
		\textsc{Institut Montpellierain Alexander Grothendieck, 499-554 Rue du Truel, 34090 Montpellier, France.}\par\nopagebreak
		\textit{E-mail address}: \texttt{anthony.genevois@umontpellier.fr}
\medskip

		\textsc{D\'epartement de Math\'ematiques B\^atiment 307, Facult\'e des Sciences d'Orsay, Universit\'e Paris Saclay, F-91405 Orsay Cedex, France.}\par\nopagebreak
		\textit{E-mail address}: \texttt{anne.lonjou@u-psud.fr}
\medskip

		\textsc{EPFL, SB MATH, Station 8, CH-1015 Lausanne, Switzerland.}\par\nopagebreak
\textit{E-mail address}: \texttt{christian.urech@epfl.ch}
\medskip

}}

\title{Asymptotically rigid mapping class groups I: Finiteness properties of braided Thompson's and Houghton's groups}
\date{\today}
\author{Anthony Genevois, Anne Lonjou, and Christian Urech}


\begin{document}

\maketitle

\begin{abstract}
This article is dedicated to the study of asymptotically rigid mapping class groups of infinitely-punctured surfaces obtained by thickening planar trees. Such groups include the braided Ptolemy-Thompson groups $T^\sharp,T^\ast$ introduced by Funar and Kapoudjian, and the braided Houghton groups $\mathrm{br}H_n$ introduced by Degenhardt. We present an elementary construction of a contractible cube complex, on which these groups act with cube-stabilisers isomorphic to finite extensions of braid groups. As an application, we prove Funar-Kapoudjian's and Degenhardt's conjectures by showing that $T^\sharp,T^\ast$ are of type $F_\infty$ and that $\mathrm{br}H_n$ is of type $F_{n-1}$ but not of type $F_n$. 
\end{abstract}

\tableofcontents

\section{Introduction}

\noindent
The groups $F$, $T$, and $V$ introduced by R. Thompson in the 1960s have a particular place in the history of group theory. First, $T$ and $V$ are the first examples of infinite finitely presented simple groups, and $F$ is the first example of a torsion-free group of type $F_\infty$ that is not of type $F$. But, since then, many groups have been constructed by varying the constructions of $F$, $T$, and $V$; see for instance \cite{HigmanV, Stein, MR1396957, Rover, Nekrashevych, BrinnV, FunarUniversal, DehornoybrV, BrinbrV, Funar-Kapoudjian, Sim, Monod, QV, MR4009393}. Although all these groups turn out to share similar properties, axiomatising the family of ``Thompson-like groups'' seems difficult; see \cite{Thumann, Witzel} for attempts in this direction. Nowadays, the investigation of Thompson-like groups is a subject on its own. Recent successes include the construction of new examples of non-amenable groups without non-abelian free subgroups \cite{Monod, LM} and the construction of simple groups distinguished by finiteness properties \cite{MR3910073, TwistedThompson}. 

\medskip \noindent
In this article, we are mainly interested in braided versions of Thompson's groups and their finiteness properties. Recall that a group $G$ is \emph{of type $F_n$} if it admits a \emph{classifying space} (i.e. an aspherical CW complex with $G$ as its fundamental group) that contains only finitely many cells in its $n$-skeleton. A group is \emph{of type $F_\infty$} if it is of type $F_n$ for every $n \geq 1$. Notice that groups of type $F_1$ coincide with finitely generated groups, and that groups of type $F_2$ coincide with finitely presented groups. Being of type $F_n$ for $n \geq 3$ is usually thought of as a higher dimensional analogue of these properties. Because the homotopy type of an aspherical CW complex depends only on its fundamental group, one can associate topological invariants to a group from a classifying space, such as (co)homology groups. Then, being of type $F_n$ assures that such invariants in dimension $\leq n$ are finitely generated. One can expect, next, to construct an explicit classifying space and to compute these invariants. Interestingly, the property of being of type $F_n$ can also be characterised from coarse geometry. In some sense, a finitely generated group is of type $F_n$ if and only if it is \emph{coarsely $(n-1)$-connected}. See \cite[Section~9.2]{MR3753580} for more details. In particular, being of type $F_n$ is a quasi-isometric invariant. 

\medskip \noindent
An interesting question is to determine to what extend a braided version of a Thompson group satisfies the same finiteness properties as the corresponding Thompson group. A positive answer in this direction can be found in \cite{MR3545879}, where the authors prove that the braided version of $V$ introduced in \cite{DehornoybrV, BrinbrV} and the braided version of $F$ introduced in \cite{MR2384840} are of type $F_\infty$. In \cite{Funar-Kapoudjian}, two braided versions of $T$, denoted by $T^\sharp$ and $T^\ast$, are introduced and proved to be finitely presented. The authors next proved that they are of type $F_3$ in \cite{MR2803858}, and they conjectured that they are of type $F_\infty$. This conjecture was the motivation of the present work.

\medskip \noindent
Our framework, largely inspired by \cite{Funar-Kapoudjian}, is the following. Fix a locally finite tree $A$ embedded into the plane in such a way that its vertex-set is discrete. The \emph{arboreal surface} $\mathscr{S}(A)$ is the oriented planar surface with boundary obtained by thickening $A$ in the plane. We denote by $\mathscr{S}^\sharp(A)$ the punctured arboreal surface obtained from $\mathscr{S}(A)$ by adding a puncture for each vertex of the tree. Now we fix a \emph{rigid structure} on $\mathscr{S}^\sharp(A)$, i.e. a decomposition into \emph{polygons} by means of a family of pairwise non-intersecting arcs whose endpoints are on the boundary of $\mathscr{S}(A)$ such that each polygon contains exactly one vertex of the underlying tree in its interior and such that each arc crosses once and transversely a unique edge of the tree. See for instance Figure~\ref{Dsharp}. We are interested in a specific subgroup of the big mapping class group of $\mathscr{S}^\sharp(A)$, corresponding to the (isotopy classes of the) homeomorphisms that, loosely speaking, preserve the rigid structure ``almost everywhere''.

\begin{definition}
A homeomorphism of $\mathscr{S}^\sharp(A)$ is \emph{asymptotically rigid} if it sends all but finitely many polygons of the rigid structure to polygons. We denote by $\amod (A)$ of isotopy classes of orientation-preserving asymptotically rigid homeomorphisms of $\mathscr{S}^\sharp(A)$.
\end{definition}

\noindent
As an example, if $A_n$ denotes the $(n+1)$-regular tree ($n \geq 2$), then $\amod(A_n)$ provides a braided version of Thompson's group $T_n$, referred to as the \emph{braided Ptolemy-Thompson group} and denoted by $\mathrm{br}T_n$. For instance, $\mathrm{br}T_2$ coincides with the group $T^\sharp$ introduced in \cite{Funar-Kapoudjian}. More generally, if $A_{n,m}$ denotes the tree with one vertex of valence $m$ while all the other vertices have valence $n+1$ ($n \geq 2$), we refer to $\mathrm{br}T_{n,m}:= \amod(A_{n,m})$ as the \emph{braided Higman-Thompson group}. So, for $m=n+1$, we recover the braided Ptolemy-Thompson group $\mathrm{br}T_n$. Interestingly, the group $T^\ast$ introduced in \cite{Funar-Kapoudjian} turns out to coincide with $\mathrm{br}T_{2,4}$; see Section \ref{sec:arboreal} for more details. As another example, if $R_n$ denotes the union of $n$ infinite rays whose origins are all identified, then the finite-index subgroup $\mathrm{br}H_n$ of elements in $\amod(R_n)$ that preserves the ends of the surface $\mathscr{S}^\sharp(R_n)$ coincides with the \emph{braided Houghton group} introduced in \cite{Degenhardt}, as observed in \cite{FunarHoughton}. Finally, let us mention that a braided version of the lamplighter group $\mathbb{Z}\wr \mathbb{Z}$ can also be constructed, but this example will be treated in more details in a forthcoming work \cite{Next}; see also Section 2.

\medskip \noindent
The major contribution of our work is the introduction of a geometric point of view on asymptotically mapping class groups. Given a simplicial tree $A$, we make $\amod(A)$ act on the cube complex $\mathscr{C}(A)$ constructed as follows. The vertices are classes of marked subsurfaces $[\Sigma, \varphi]$ where $\Sigma$ is a non-empty connected union of finitely many polygons of the rigid structure and $\varphi$ is an asymptotically rigid homeomorphism of $\mathscr{S}^\sharp(A)$ (for a more precise definition see Section \ref{section:Construction}). For any $H_1, \ldots, H_k$ pairwise distinct adjacent polygons of $\Sigma$, there is a $k$-cube issued from $[\Sigma, \varphi]$ spanned by $$\left\{ \left[  \Sigma \cup \bigcup\limits_{i \in I} H_i , \varphi \right] \mid I \subset \{1, \ldots, k\} \right\}.$$
This simple cube complex turns out to be tightly related to $\amod(A)$. As a first application, we show that every finite subgroup in $\amod(A)$ fixes a vertex in $\mathscr{C}(A)$, which allows us to prove that finite subgroups in $\amod(A)$ are all cyclic with specific orders. See Theorem \ref{thm:FiniteOrders}. As a particular case, we characterise all the orders of finite-order elements in braided Higman-Thompson groups:

\begin{restatable*}{thm}{torsionthm}\label{prop:torsion}
Let $n \geq 2$ and $m \geq 1$ be two integers. The braided Higman-Thompson group $\mathrm{br}T_{n,m}$ contains an element of order $l$ if and only if 
$$\left\{ \begin{array}{ll} \text{$l$ divides $m$ or $m-n+1$} & \text{if $m \neq n+1$} \\ \text{$l=2$ or $l$ divides $n+1$} & \text{if $m=n+1$}  \end{array} \right..$$
\end{restatable*}

\noindent
This result has interesting consequences on the isomorphism problem. For instance, we deduce that the braided Ptolemy-Thompson groups $\mathrm{br}T_n$ and $\mathrm{br}T_m$ are isomorphic if and only if $n=m$. We also recover the fact, proved in \cite{Funar-Kapoudjian}, that $T^\sharp$ and $T^\ast$ are not isomorphic. However, not all braided Higman-Thompson groups can be distinguished by their finite-order elements: notice that, for every $n \geq 2$, $\mathrm{br}T_{n,n-1}$ contains an element of every possible order.

\medskip \noindent
But our main application concerns finiteness properties. Such applications are motivated by the following observation (which is a consequence of Theorem \ref{thm:contractible} and Lemma \ref{lem:VertexStab}):

\begin{prop}
For every simplicial tree $A$, $\mathscr{C}(A)$ is a contractible cube complex, on which $\amod(A)$ acts with cube-stabilisers isomorphic to finite extensions of braid groups.
\end{prop}

\noindent
As a consequence, by looking at the action of the groups $\amod(A)$ on the cube complexes $\mathscr{C}(A)$, we are in a good position to prove finiteness properties among asymptotically rigid mapping class groups. The two main theorems of this article go in this direction:

\begin{restatable*}{thm}{ThompsonThm}\label{PtolemyConnected}
For all $n \geq 2$ and $m \geq 1$, the braided Higman-Thompson group $\mathrm{br}T_{n,m}$ is of type $F_\infty$.
\end{restatable*}

\begin{restatable*}{thm}{HoughtonThm}\label{thm:brHfiniteness}
For every $n \geq 1$, the braided Houghton group $\mathrm{br}H_n$ is of type $FP_{n-1}$ but not of type $FP_n$. Moreover, if $n \geq 3$, $\mathrm{br}H_n$ is finitely presented.
\end{restatable*}

\noindent
The property $FP_n$ is a cohomological analogue of the property $F_n$. A group of type $F_n$ is automatically of type $FP_n$, and the converse holds for finitely presented groups (see for instance \cite[Section 8.7]{MR1324339}). Consequently, Theorem \ref{thm:brHfiniteness} shows that, for every $n \geq 1$, the braided Houghton group $\mathrm{br}H_n$ is of type $F_{n-1}$ but not of type $F_n$.

\medskip \noindent
For $(n,m)=(2,1)$ and $(n,m)=(2,4)$, Theorem \ref{PtolemyConnected} proves Funar and Kapoudjian's conjecture according to which $T^\sharp$ and $T^\ast$ are of type $F_\infty$. Theorem \ref{thm:brHfiniteness} was conjectured in \cite{Degenhardt} and verified for $n \leq 3$. A strategy was suggested in \cite{BuxHoughton} for the general case, but our proof of the conjecture uses a different approach. Theorems \ref{PtolemyConnected} and \ref{thm:brHfiniteness} show that the braided versions of Thompson's and Houghton's groups verify the same finiteness properties as the unbraided versions.

\medskip \noindent
We emphasize that our proofs of Theorems \ref{PtolemyConnected} and \ref{thm:brHfiniteness} are constructive, providing explicit highly connected complexes on which asymptotically rigid mapping class groups act. As an application, in a forthcoming article \cite{NextNext}, we compute explicit presentations of the braided Higman-Thompson groups. In particular, this allows us to compute the abelianisations of these groups, providing other algebraic invariants in addition to Theorem \ref{prop:torsion}. We also emphasize that the techniques developped in this article go beyond the arboreal surfaces we study here. For instance, we expect that our arguments can be adapted to the universal mapping class groups constructed in \cite{MR2105950, AramayonaFunar}.

\paragraph{A few words about the proofs.} 
If we allow $n=1$ in our notation $A_{n,m}$, then, when $n \geq 2$, $\amod(A_{n,m})$ coincides with the braided Higman-Thompson group $\mathrm{br}T_{n,m}$; and, when $n=1$, $\amod(A_{n,m})$ contains the braided Houghton group $\mathrm{br}H_m$ as a finite-index subgroup. This notation allows us to work with these two families of groups simultaneously. The cube complex $\mathscr{C}(A_{n,m})$, on which $\amod(A_{n,m})$ acts, is naturally endowed with the \emph{height function}
$$h : \text{vertex $[\Sigma, \varphi]$} \mapsto \text{number of punctures in $\Sigma$};$$
and standard arguments from Morse theory allow us to deduce finiteness properties of $\amod(A_{n,m})$ from a careful analysis of the descending links in $\mathscr{C}(A_{n,m})$. However, this strategy may fail because of vertices of arbitrarily large height with non-simply connected descending links. We first need to extract from $\mathscr{C}(A_{n,m})$ a suitable subcomplex $\mathscr{SC}(A_{n,m})$, referred to as the \emph{spine}. It corresponds to the subcomplex spanned by the vertices of the form $[\Sigma,\varphi]$ where $\Sigma$ contains the vertex of valence $m$ in $A_{n,m}$, and it turns out to be homotopy equivalent to the original complex. In particular, it is also contractible. 

\medskip \noindent
Interestingly, the descending links in $\mathscr{SC}(A_{n,m})$ can be described as complexes of arcs $\mathfrak{C}(p,q,r)$ for some $p,q \geq 1$ and $r \geq 0$.

\begin{definition}
Let $p, q\geq1$ and $r \geq 0$ be three integers. Fix a disc $\mathbb{D}$ with $p$ punctures in its interior and $q$ marked points on its boundary. Let $\{m_i \mid i \in \mathbb{Z}_q \}$ denote these marked points, ordered cyclically. From now on, an arc in $\mathbb{D}$ refers to an arc that starts from a marked point and that ends at a puncture. Two arcs that start from the marked points $m_i,m_j$ are \emph{$r$-separated} if they are disjoint and if the distance between $i$ and $j$ in $\mathbb{Z}_q$ is $>r$ (where $\mathbb{Z}_q$ is metrically thought of as the cycle $\mathrm{Cayl}(\mathbb{Z}_q,\{1\})$). Notice that being $0$-separated amounts to being disjoint. We define $\mathfrak{C}(p,q,r)$ as the simplicial complex whose vertices are the isotopy classes of arcs and whose simplices are collections of arcs that are pairwise $r$-separated (up to isotopy). 
\end{definition}

\noindent
Then, proving Theorems \ref{PtolemyConnected} and \ref{thm:brHfiniteness} essentially amounts to showing that, for a fixed $r \geq 0$, our arc complex $\mathfrak{C}(p,q,r)$ becomes more and more connected as $p$ and $q$ increase. More precisely, we prove that:

\begin{thm}\label{Intro:ArcComplex}
Let $p \geq 2$, $q \geq 1$ and $r \geq 0$ be three integers. The complex $\mathfrak{C}(p,q,r)$ is $\left( \left\lfloor \frac{1}{3} \left( p+ \left\lfloor \frac{q}{r+1} \right\rfloor \right) \right\rfloor -2 \right)$-connected. Moreover, if $r=0$, then $\mathfrak{C}(p,q,r)$ is homotopy equivalent to a bouquet of infinitely many $(q-1)$-spheres as soon as $p \geq 2q$.
\end{thm}

\noindent
Our theorem, proved more generally for arbitrary surfaces and not only for discs, is a direct consequence of Propositions \ref{prop:SimConnected} and \ref{prop:BouquetSpheres} below. Our argument follows closely the lines of the proof of \cite[Theorem~3.10]{MR3545879}, which is inspired by an argument from \cite{MR1123262}. 

\medskip \noindent
The topology of complexes of arcs has been widely studied in the literature, especially because of its connection with homology stability. See for instance the foundational article \cite{MR786348}. However, it is usually assumed that a collection of arcs that have disjoint interiors but that start from the same marked point on the boundary spans a simplex, which is not the case in $\mathfrak{C}(p,q,r)$. Even worse, in order to span a simplex, we require the marked points to be sufficiently far away from each other. So $\mathfrak{C}(p,q,r)$ is significantly different from the complexes of arcs usually studied.

\begin{remark}
After the completion of this work, K.-U. Bux informed us that he also had a proof of Theorem \ref{thm:brHfiniteness}. His proof is announced in \cite[Remark 13.10.2]{BuxBraided} and some details are given in \cite{BuxArcComplexes}. Notice that the complex $\mathcal{A}(m,n,0)$ defined in \cite{BuxArcComplexes} coincides with our $\mathfrak{C}(m,n,0)$, so the families of arc complexes considered in \cite{BuxArcComplexes} and in our paper overlap. However, no one is included in the other. In particular, $\mathfrak{C}(m,n,r)$ for $r \geq 1$ does not belong to the framework of \cite{BuxArcComplexes}. 
\end{remark}

\paragraph{Organisation of the article.} Section~\ref{sec:arboreal} is dedicated to basic definitions and examples about asymptotically rigid mapping class groups. Our cube complex, the main object of the article, is constructed in Section~\ref{section:CC}, where we prove that it is contractible. A quick discussion related to its curvature is included in Subsection~\ref{section:curvature}. As a first application, we classify finite subgroups in asymptotically rigid mapping class groups in Section~\ref{section:torsion}. The core of the article is Section~\ref{section:TypeF} where we prove Theorems~\ref{PtolemyConnected} and~\ref{thm:brHfiniteness}. The spine is introduced and studied in Subsection~\ref{section:spine} and the descending links of its vertices are described as arc complexes in Subsection~\ref{section:Links}. Finally, we study their homotopy (Theorem~\ref{Intro:ArcComplex}) and we prove Theorems~\ref{PtolemyConnected} and~\ref{thm:brHfiniteness} in Subsections~\ref{section:HigmanThompson} and~\ref{section:Houghton}.

\paragraph{Acknowledgments.} The authors thank the University of Basel and the EPFL for their hospitality during parts of this project. A. G. was supported by a public grant as part of the Fondation Math\'ematique Jacques Hadamard. A. L. acknowledges support from the Swiss National Science
Foundation Grant ``Birational transformations of threefolds''
$200020\mathunderscore178807$.

\section{Arboreal surfaces and asymptotically rigid mapping class groups}\label{sec:arboreal}

\noindent
Let us recall from the introduction the general framework of this article. 

\begin{figure}
\begin{center}
\includegraphics[scale=0.3]{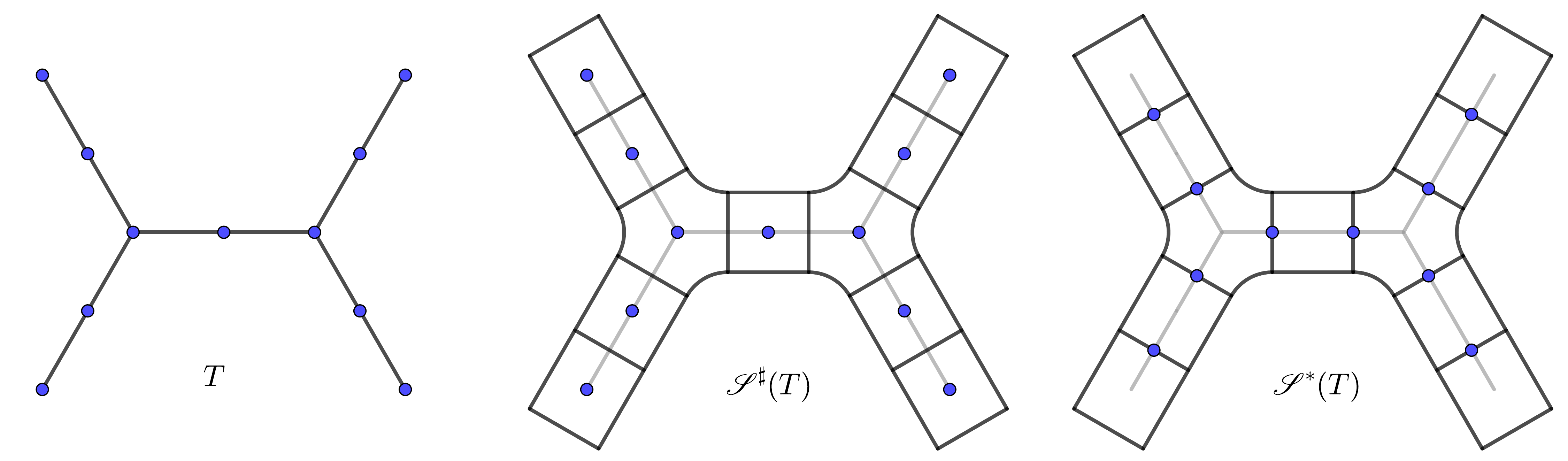}
\caption{Surfaces with rigid structures associated to a simplicial tree.}
\label{Dsharp}
\end{center}
\end{figure}

\medskip \noindent
Fix a locally finite tree $A$ embedded into the plane in such a way that its vertex-set is discrete. The \emph{arboreal surface} $\mathscr{S}(A)$ is the oriented planar surface with boundary obtained by thickening $A$ in the plane. We denote by $\mathscr{S}^\sharp(A)$ the punctured arboreal surface obtained from $\mathscr{S}(A)$ by adding a puncture for each vertex of the tree. Following \cite{Funar-Kapoudjian}, we fix a \emph{rigid structure} on $\mathscr{S}^\sharp(A)$, i.e. a decomposition into \emph{polygons} by means of a family of pairwise non-intersecting arcs whose endpoints are on the boundary of $\mathscr{S}(A)$ such that each polygon contains exactly one vertex of the underlying tree in its interior and such that each arc crosses once and transversely a unique edge of the tree. See for instance Figure~\ref{Dsharp}.

\medskip \noindent
A subsurface of $\mathscr{S}^\sharp(A)$ is \emph{admissible} if it is a non-empty connected finite union of polygons belonging to the rigid structure. A homeomorphism $\varphi : \mathscr{S}^\sharp(A) \to \mathscr{S}^\sharp(A)$ is \emph{asymptotically rigid} if the following conditions are satisfied:
\begin{itemize}
	\item there exists an admissible subsurface $\Sigma \subset \mathscr{S}^\sharp(A)$ such that $\varphi(\Sigma)$ is also admissible;
	\item the homeomorphism $\varphi$ is \emph{rigid outside $\Sigma$}, i.e. the restriction \[\varphi : \mathscr{S}^\sharp(A) \backslash \Sigma \to \mathscr{S}^\sharp(A) \backslash \varphi( \Sigma)\]
	 respects the rigid structure, mapping polygons to polygons. Such a surface $\Sigma$ is called a \emph{support} for $\varphi$.
\end{itemize}
We denote by $\amod (A)$ the group of isotopy classes of orientation-preserving asymptotically rigid homeomorphisms of $\mathscr{S}^\sharp(A)$. We emphasize that isotopies have to fix each puncture. 

\medskip \noindent
In the sequel, we refer to the \emph{frontier} $\mathrm{Fr}(\Sigma)$ of an admissible subsurface $\Sigma$ as the union of the arcs defining the rigid structure that are contained in the boundary. Also, a polygon is called \emph{adjacent} to $\Sigma$ if it is not contained in $\Sigma$ but shares an arc with the frontier of $\Sigma$.

\medskip \noindent
Any asymptotically rigid homeomorphism $\mathscr{S}^\sharp(A) \to \mathscr{S}^\sharp(A)$ induces a \emph{quasi-auto\-morph\-ism} $A \to A$, i.e. a bijection defined on the vertices of $A$ that preserves adjacency and non-adjacency for all but finitely many pairs of vertices. Let $\mathrm{QAut}(A)$ denote the group of quasi-automorphisms of $A$. The induced morphism $\amod(A) \to \mathrm{QAut}(A)$ is not surjective in general, and we denote its image by $\mathfrak{mod}_a(A)$; its kernel corresponds to the homeomorphisms $\mathscr{S}^\sharp(A) \to \mathscr{S}^\sharp(A)$ that are the identity outside an admissible subsurface and that fix each puncture. In other words, we have a short exact sequence
$$1 \to PB_\infty \to \amod(A) \to \mathfrak{mod}_a(A) \to 1$$
where $PB_\infty$ denotes the limit of pure mapping class groups $\bigcup\limits_{\text{$\Sigma$ admissible}} \mathrm{PMod}(\Sigma)$. Notice that $PB_\infty$ is isomorphic to the group of compactly supported pure braids with infinitely many strands. We refer to the previous exact sequence as the \emph{arboreal short exact sequence} satisfied by $\amod(A)$.

\medskip \noindent
Another interesting exact sequence comes from the \emph{forgetful map} $\amod(A) \to \Mod(\mathscr{S}(A))$, where we forget the punctures. We denote by $\amod_f(A)$ the image of $\amod(A)$ in the mapping class group $\Mod(\mathscr{S}(A))$. The kernel of the forgetful map corresponds to the homeomorphisms $\mathscr{S}^\sharp(A) \to \mathscr{S}^\sharp(A)$ that are the identity outside an admissible subsurface. In other words, we have a short exact sequence
$$1 \to B_\infty \to \amod(A) \to \amod_f(A) \to 1$$
where $B_\infty$ denotes the limit $\bigcup\limits_{\text{$\Sigma$ admissible}} \mathrm{Mod}(\Sigma)$. Notice that $B_\infty$ is isomorphic to the group of compactly supported braids with infinitely many strands. We refer to the previous exact sequence as the \emph{forgetful short exact sequence} satisfied by $\amod(A)$.

\medskip \noindent
These two short exact sequences motivate the idea that $\amod(A)$ is a braided version of $\amod_a(A)$ and $\amod_f(A)$. They may also be useful to compute explicit presentations; see Example \ref{ex:brH2}. However, as $PB_\infty$ and $B_\infty$ are not even finitely generated, such presentations are infinite, and in general it is not clear if finite presentations can be extracted from them, and if so, how. 

\medskip \noindent
Depending on the choice of the tree $A$, we obtain various groups with a rich structure. Let us present some particular cases of groups that interest us most.

\paragraph{Braided Higman-Thompson groups.}
For integers $n\geq 2$ and $m\geq 1$, let $A_{n,m}$ be the tree with one vertex of valence $m$ while all the other vertices have valence $n+1$. We then call the group $\amod(A_{n,m})$ the \emph{braided Higman-Thompson group} and denote it by $\mathrm{br}T_{n,m}$. The terminology is justified by the forgetful short exact sequence
$$1 \to B_\infty \to \mathrm{br}T_{n,m} \to T_{n,m} \to 1,$$
where $T_{n,m}:= \amod_f(A_{n,m})$ coincides with the Higman-Thompson group of type $(n,m)$ introduced in \cite{Brown} by analogy with \cite{HigmanV}. (We refer to \cite[Section 2.3]{Survey} for a justification of this identification for $n=2$ and $m=1$. The general case follows similarly.)
Note that $A_n:=A_{n,n+1}$ is the $(n+1)$-regular tree. We refer to $\mathrm{br}T_n:= \mathrm{br}T_{n,n+1}$ as the \emph{braided Ptolemy-Thompson group}. For $n=2$, we recover the group $T^\sharp$ introduced in \cite{Funar-Kapoudjian}.

\begin{ex}
Notice that many isometries of $A_n$ induce (asymptotically) rigid homeomorphisms of $\mathscr{S}^\sharp(A_n)$. Moreover, any two distinct such isometries induce non-isotopic homeomorphisms, so many subgroups of $\mathrm{Aut}(A_n)$ turn out to define subgroups of $\mathrm{br}T_n$. They include non-abelian free subgroups (which do not come from braid subgroups) and finite cyclic subgroups. More precisely, if $H$ is a polygon of the rigid structure, then the homeomorphism $r_H$ that shifts cyclically the arcs of the frontier of $H$ defines a(n asymptotically) rigid homeomorphism of $\mathscr{S}^\sharp(A_n)$ of finite order.
\end{ex}

\begin{figure}
\begin{center}
\includegraphics[scale=0.4]{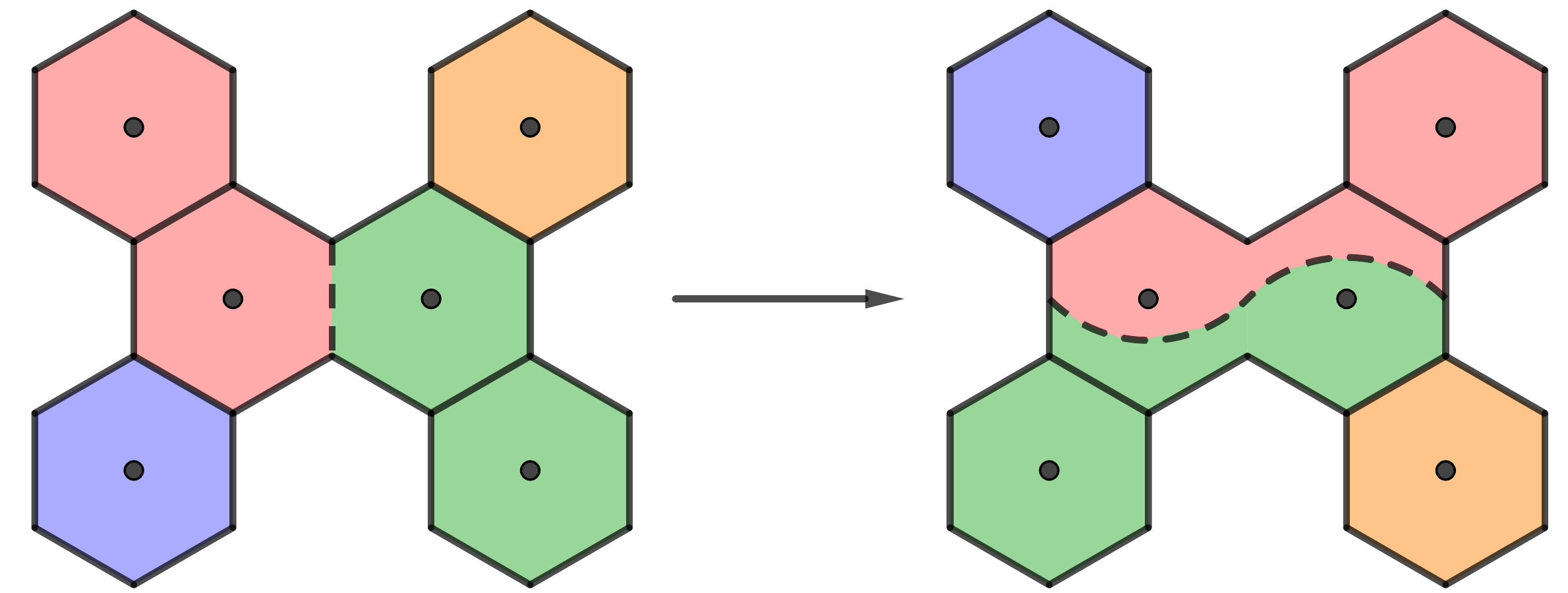}
\caption{An asymptotically rigid homeomorphism.}
\label{rotation}
\end{center}
\end{figure}

\begin{ex}\label{ex_rotations}
More generally, if $\Sigma$ is any admissible subsurface containing $k$ punctures, then its frontier consists of $k(n-1)+2$ arcs. Consequently, the complement of $\Sigma$ in  $\mathscr{S}^\sharp(A_n)$ consists of $k(n-1)+2$ pairwise homeomorphic arboreal surfaces. We denote by $r_\Sigma$ the asymptotically rigid homeomorphism that cyclically clockwise shifts the arcs of the frontier of $\Sigma$ (and hence the homeomorphic arboreal surfaces) and whose restriction to a disk in $\Sigma$ containing all the punctures is the identity. Figure \ref{rotation} illustrates the case $k=2$.
\end{ex}

\vspace{-0.3cm}

\paragraph{Braided Houghton groups.}
Fix an $n \geq 1$ and let $R_n$ denote the union of $n$ infinite rays whose initial vertices are identified. We refer to the index-$n$ subgroup of $\amod(R_n)$ that stabilises each end of the surface $\mathscr{S}(R_n)$ as the \emph{braided Houghton group} $\mathrm{br}H_n$. As explained below, $\mathrm{br}H_n$ coincides with the groups introduced in \cite{Degenhardt} and \cite{FunarHoughton}. The terminology is justified by the exact sequence
$$1 \to PB_\infty \to \mathrm{br}H_n \to H_n \to 1$$
given by the arboreal short exact sequence satisfied by $\amod(R_n)$, where $H_n$ (defined as the finite-index subgroup in $\mathrm{QAut}(R_n)$ that stabilises each end of $R_n$) denotes the Houghton group as introduced in \cite{Houghton}. 

\begin{ex}\label{ex:brH2}
The braided Houghton group $\mathrm{br}H_2$ satisfies the short exact sequence
$$1 \to B_\infty \to \mathrm{br}H_2 \to \mathbb{Z} \to 1$$
where the cyclic group corresponds to the subgroup in $\amod_f(R_2)$ that fixes the two ends of $R_2$. This exact sequence splits and provides the decomposition $\mathrm{br}H_2 = B_\infty \rtimes \langle t \rangle$, where $t$ denotes the homeomorphism coming from the translation acting on $R_2$. It follows that $\mathrm{br}H_2$ admits
$$\left\langle t, \tau_i \ (i \in \mathbb{Z}) \left| \begin{array}{l} [\tau_i, \tau_j]=1, \ i,j \in \mathbb{Z}, \ |i-j| \geq 2 \\ \tau_i \tau_{i+1} \tau_i = \tau_{i+1} \tau_i \tau_{i+1}, \ i \in \mathbb{Z} \\ \tau_i^t=\tau_{i+1}, \ i \in \mathbb{Z} \end{array} \right\rangle \right.$$
as a presentation, where $\tau_i$ corresponds to the twist between the punctures $i$ and $i+1$, and with the notation $a^b:= bab^{-1}$. Setting $\tau:= \tau_0$ and using the relation $\tau_i = t^i \tau t^{-i}$ for every $i \in \mathbb{Z}$, the presentation can be simplified as:
$$\left\langle t, \tau \mid \tau \tau^t \tau = \tau^t \tau \tau^t, \ \left[ \tau, \tau^{t^n} \right]=1 \ (n \in \mathbb{Z}, |n| \geq 2) \right\rangle.$$
(As proved in \cite{Degenhardt}, and reproved by Theorem~\ref{thm:brHfiniteness}, $\mathrm{br}H_2$ is not finitely presented, so a finite presentation cannot be extracted from the previous presentation.)
\end{ex}

\paragraph{Braided lamplighter group.}
Let $A$ denote the tree obtained from an horizontal bi-infinite line by gluing an infinite descending vertical ray to each of its vertices. We refer to $\mathrm{br}\mathcal{L}:= \amod (A)$ as the \emph{braided lamplighter group}. The terminology is justified by the short exact sequence
$$1 \to B_\infty \to \mathrm{br}\mathcal{L} \to \mathcal{L} \to 1$$
given by the forgetful short exact sequence satisfied by $\amod(A)$, where $\mathcal{L}$ denotes the \emph{lamplighter group}, defined as the wreath product $\mathbb{Z} \wr\mathbb{Z}:= \left( \bigoplus\limits_\mathbb{Z} \mathbb{Z} \right) \rtimes \mathbb{Z}$ (where $\mathbb{Z}$ acts on the direct sum by shifting the coordinates). In a forthcoming article \cite{Next}, it will be proved that the braided lamplighter group $\mathrm{br}\mathcal{L}$, like its unbraided version $\mathcal{L}$, is finitely generated but not finitely presented.

\paragraph{Variations of the definition.} Following \cite{Funar-Kapoudjian}, given a simplicial tree $A$ we can define a different surface with rigid structure $\mathscr{S}^\ast(A)$ from $\mathscr{S}(A)$: we add a puncture for each edge of the tree, and we decompose the surface into polygons by means of a family of pairwise non-intersecting arcs whose endpoints are on the boundary of $\mathscr{S}(A)$ such that each arc contains a puncture and crosses the underlying tree once and transversely. See Figure \ref{Dsharp}. Mimicking the definition of $\amod(A)$, one defines a new group $\amod^\ast(A)$. As observed in \cite{Funar-Kapoudjian}, the groups $\amod(A)$ and $\amod^\ast(A)$ may not be isomorphic (see Remark~\ref{rem:funar_kap_non-iso} below). However, as justified by the next observation, there is no loss of generality in studying only the groups $\amod(A)$.

\begin{lemma}\label{lem:Iso}
Let $A$ be a simplicial tree. If $A'$ is a tree obtained from $A$ by collapsing an edge, then $\amod^\ast(A)$ is isomorphic to $\amod(A')$.
\end{lemma}

\begin{proof}
Fix a vertex $u \in A$. Given a polygon $H$ of the rigid structure of $\mathscr{S}^\sharp(A)$ that doesn't contain $u$, let $\alpha$ denote the arc in $\mathrm{Fr}(H)$ that separates $H\backslash \alpha$  from the puncture $u$.
There is an isotopy of $\mathscr{S}(A)$ supported in $(H \backslash N(\mathrm{Fr}(H))) \cup N(\alpha)$, where $N(\alpha)$ is a small tubular neighborhood of $\alpha$ and $N(\mathrm{Fr}(H))$ is a small tubular neighborhood of $\mathrm{Fr}(H)$, that sends $\alpha$ to an arc that passes through the puncture of $H$ and that crosses the underlying tree once and transversely. Since all these isotopies have pairwise disjoint supports, we can perform on all the polygons at once and, after removing the puncture $u$, we obtain a homeomorphism $\mathscr{S}^\sharp(A) \cup \{u\} \to \mathscr{S}^\ast(A)$ that respects the rigid structures, sending polygons to polygons. (This paragraph is extracted from the proof of \cite[Proposition~2.10]{Funar-Kapoudjian}.) 

\medskip \noindent
Now, fix a neighbor $v \in A$ of $u$ and let $A'$ denote the tree obtained from $A$ by collapsing the edge between $u$ and $v$. Let $A,B$ denote the polygons of $\mathscr{S}^\sharp(A)$ containing the punctures $u,v$. Define a new rigid structure on $\mathscr{S}^\sharp(A) \cup \{u\}$ by removing the arc common to $A,B$. (In other words, we are merging $A$ and $B$ into a unique polygon.) Notice that the surface with rigid structure we obtain coincides with $\mathscr{S}^\sharp(A')$ (up to a homeomorphism that respects the rigid structure).

\medskip \noindent
Therefore, there exists a homeomorphism $\mathscr{S}^\sharp(A') \to \mathscr{S}^\ast(A)$ that sends each polygon of $\mathscr{S}^\sharp(A')$ to a polygon of $\mathscr{S}^\ast(A)$, except one that is sent to the union of two polygons. By conjugating by such a homeomorphism, one obtains an isomorphism $\amod(A') \to \amod^\ast(A)$. 
\end{proof}

\noindent
In \cite{Funar-Kapoudjian}, the authors study the groups $T^\sharp:= \amod(A_2)$ and $T^\ast:= \amod^\ast(A_2)$ and prove that they are not isomorphic by comparing their abelianisations (see also Remark \ref{rem:funar_kap_non-iso}). As a consequence of Lemma~\ref{lem:Iso}, $T^\ast$ turns out to be isomorphic to $\amod(A_{2,4})$. Thus, one may justify the difference between $T^\sharp$ and $T^\ast$ by saying that $T^\sharp$ is a braided version of the Ptolemy-Thompson group $T_2$ but that $T^\ast$ is a braided version of the Higman-Thompson group $T_{2,4}$ (as defined in \cite{Brown}, following \cite{HigmanV}).

\begin{remark}
As a consequence of Lemma \ref{lem:Iso}, if $A_1$ and $A_2$ are two trees such that a third tree $A$ can be obtained from both $A_1$ and $A_2$ by collapsing an edge, then $\amod(A_1)$ and $\amod(A_2)$ are isomorphic. Therefore, we obtain non-trivial examples of isomorphic asymptotically rigid mapping class groups.
\end{remark}

\begin{figure}
\begin{center}
\includegraphics[scale=0.28]{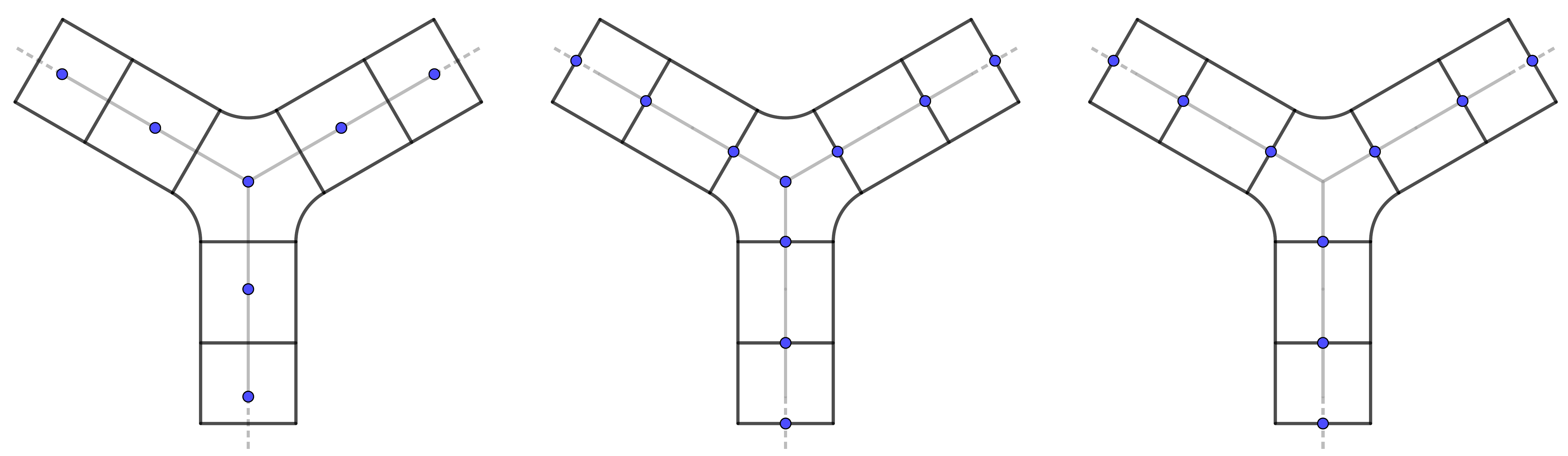}
\caption{}
\label{houghtons}
\end{center}
\end{figure}

\begin{remark}
It is worth noticing that, as a consequence of Lemma \ref{lem:Iso} (and its proof), the three braided versions of the Houghton groups associated to the surfaces with rigid structures illustrated by Figure \ref{houghtons} all coincide. We used the first surface in our definition of braided Houghton's groups; the second surface is used in \cite{FunarHoughton} in order to recover the original definition of \cite{Degenhardt}; and the third surface is used in \cite{FunarHoughton} in order to define an a priori different braided version of Houghton's groups. 
\end{remark}

\section{A contractible cube complex}\label{section:CC}

\subsection{Construction}\label{section:Construction}

\noindent
In this section, we fix a locally finite tree $A$ embedded into the plane in such a way that its vertex-set is discrete. Let $\mathscr{C}(A)$ denote the cube complex defined as follows:
\begin{itemize}
	\item A vertex of $\mathscr{C}(A)$ is a couple $(\Sigma,\varphi)$, where $\Sigma \subset \mathscr{S}^\sharp(A)$ is an admissible subsurface and $\varphi : \mathscr{S}^\sharp(A) \to \mathscr{S}^\sharp(A)$ an asymptotically rigid homeomorphism, modulo the following equivalence relation: $(\Sigma_1, \varphi_1) \sim (\Sigma_2, \varphi_2)$ if $\varphi_2^{-1} \varphi_1$ is isotopic to an asymptotically rigid homeomorphism  that maps $\Sigma_1$ to $\Sigma_2$ and that is rigid outside $\Sigma_1$. We denote by $[\Sigma,\varphi]$ the vertex of $\mathscr{C}(A)$ represented by $(\Sigma, \varphi)$.
	\item An edge of $\mathscr{C}(A)$ links $[\Sigma, \varphi]$ and $[\Sigma \cup H, \varphi]$ where $H$ is a polygon adjacent to $\Sigma$.
	\item If $[\Sigma, \varphi]$ is a vertex and if $H_1, \ldots, H_k$ are pairwise distinct adjacent polygons of $\Sigma$, the subgraph spanned by $$\left\{ \left[  \Sigma \cup \bigcup\limits_{i \in I} H_i , \varphi \right] \mid I \subset \{1, \ldots, k\} \right\}$$ is filled in with a $k$-cube.
\end{itemize}

\begin{figure}
\begin{center}
\includegraphics[trim=0 13cm 0 0,clip,scale=0.27]{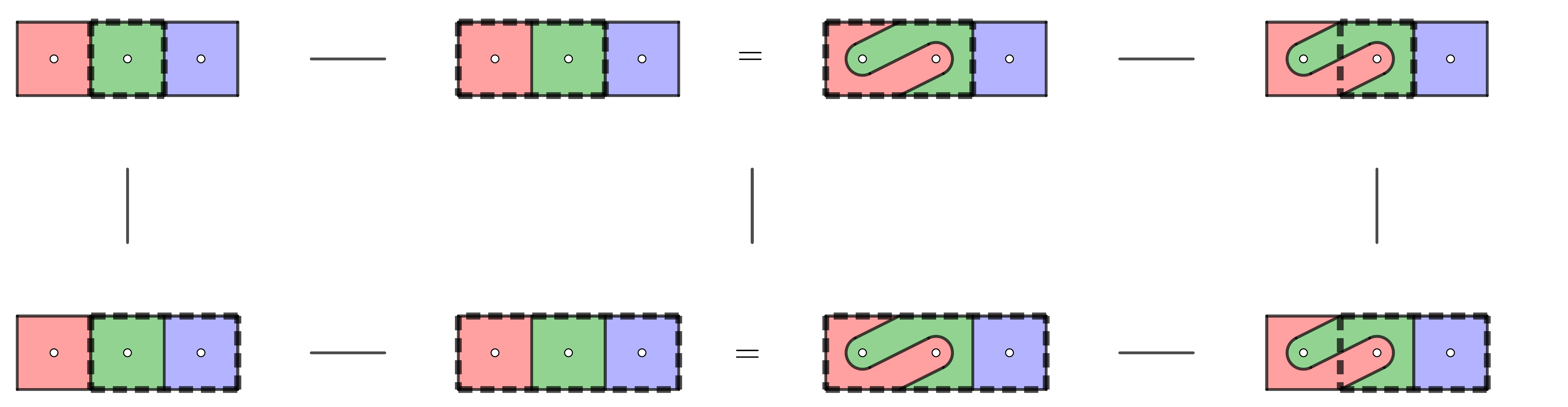}
\caption{Two adjacent squares in the cube complex $\mathscr{C}(A)$ when $A$ is a linear tree of length two. The colored $2$-cells indicate the markings. So the vertices in the left square are marked by the identity and the two right vertices are marked by a twist.}
\label{ExCC}
\end{center}
\end{figure}

\noindent
Figure \ref{ExCC} illustrates a piece of $\mathscr{C}(A)$ when $A$ is a linear tree of length two. The group $\amod(A)$ naturally acts on the cube complex $\mathscr{C}(A)$ by isometries preserving the structure of a cube complex in the following way: for every asymptotically rigid homeomorphism $g \in \amod(A)$ and every $[\Sigma, \varphi] \in \mathscr{C}(A)$, we define
$$g \cdot [\Sigma, \varphi]:= [\Sigma, g \varphi].$$
Let us note that this action is well-defined. If $g_1,g_2$ are two representatives of $g$ in $\amod(A)$ and if $(\Sigma_1,\varphi_1),(\Sigma_2, \varphi_2)$ are two representatives of the vertex $[\Sigma,\varphi]$, then $g_2^{-1}g_1$ is isotopic to the identity and $\varphi_2^{-1} \varphi_1$ is isotopic to a homeomorphism that sends $\Sigma_1$ to $\Sigma_2$ and that is rigid outside $\Sigma_1$. So $\varphi_2^{-1} g_2^{-1}g_1 \varphi_1$ is isotopic to a homeomorphism that sends $\Sigma_1$ to $\Sigma_2$ and that is rigid outside $\Sigma_1$, i.e. $[\Sigma_1,g_1 \varphi_1]= [\Sigma_2, g_2 \varphi_2]$. 

\medskip \noindent
It is worth noticing that our cube complex is naturally endowed with a \emph{Morse function}, which will be a key tool to prove the contractibility of $\mathscr{C}(A)$. It is defined as follows.

\medskip \noindent
Observe that, if $[\Sigma_1,\varphi_1]=[\Sigma_2,\varphi_2]$, then the surfaces $\Sigma_1$ and $\Sigma_2$ are homeomorphic, so they must have the same number of punctures. This allows us to define the \emph{height} of a vertex $x=[\Sigma, \varphi]$ as the number of punctures contained in $\Sigma$; we denote it by $h(x)$. Notice that, by construction of $\mathscr{C}(A)$, if $x$ and $y$ are two adjacent vertices then $h(y)=h(x) \pm 1$. Hence, the edges of $\mathscr{C}(A)$ are naturally oriented by the height function (from small to large height). Also, notice that the action of $\amod(A)$ preserves the height function.

\medskip \noindent
Let us record a few elementary observations about the cube complex $\mathscr{C}(A)$.

\begin{claim}
If $A$ contains at least two vertices, then $\mathscr{C}(A)$ is not locally compact.
\end{claim}

\noindent
Consider a vertex $[\Sigma, \id]$ of height $2$, $H$ one of the two polygons of $\Sigma$, and let $\tau\in\amod(A)$ correspond to the homeomorphism that is the identity outside $\Sigma$, but twists the two punctures inside $\Sigma$. Then $[\Sigma, \id]=[\Sigma, \tau^n]$ for all $n\in\Z$, but $[\Sigma\setminus H, \tau^n]\neq [\Sigma\setminus H, \tau^m]$ for all $m\neq n$. Thus we obtain infinitely many edges descending from $[\Sigma, \id]$ to vertices of smaller height.

\begin{claim}
The cube complex $\mathscr{C}(A)$ is locally finite-dimensional.
\end{claim}

\noindent
Fix a vertex $x:= [\Lambda,\psi] \in \mathscr{C}(A)$ and assume that it belongs to an $n$-cube for some $n \geq 1$. By construction, there exist an admissible subsurface $\Sigma$, polygons $H_1, \ldots, H_n$ and a homeomorphism $\varphi$ such that our cube is spanned by the vertices
$$\left\{ \left[  \Sigma \cup \bigcup\limits_{i \in I} H_i , \varphi \right] \mid I \subset \{1, \ldots, n\} \right\}.$$
As a consequence, $x=[\Sigma \cup H_{i_1} \cup \cdots \cup H_{i_r},\varphi]$ for some $1 \leq i_1, \ldots, i_r \leq n$ and $0\leq r \leq n$. Because we also have $x=[\Lambda,\psi]$, there must exist a homeomorphism of $\mathscr{S}^\sharp(A)$ sending $\Lambda$ to $\Sigma \cup H_{i_1} \cup \cdots \cup H_{i_r}$; as a consequence, the frontiers of these two surfaces have the same number of connected components. But the dimension $n$ of our cube is necessarily bounded above by the number of components of the frontier of $\Sigma$, which is at most the number of components of the frontier of $\Sigma \cup H_{i_1} \cup \cdots \cup H_{i_r}$. We conclude that any cube in $\mathscr{C}(A)$ containing $x$ has dimension at most the number of components in the frontier of~$\Lambda$.

\subsection{Contractibility}

\noindent
In this section, we fix a locally finite tree $A$ embedded into the plane in such a way that its vertex-set is discrete and we show that the complex $\mathscr{C}(A)$ constructed in Section \ref{section:Construction} is contractible:

\begin{thm}\label{thm:contractible}
The cube complex $\mathscr{C}(A)$ is contractible.
\end{thm}

\noindent
Our first preliminary lemma is an elementary statement, which will be fundamental in the sequel.

\begin{lemma}\label{lem:Choice}
Let $x,y$ two adjacent vertices in $\mathscr{C}(A)$ such that $h(y)>h(x)$. For every representative $(\Sigma, \varphi)$ of $x$, there exists a polygon $H$ adjacent to $\Sigma$ such that $y=[\Sigma \cup H, \varphi]$. 
\end{lemma}

\begin{proof}
There exist an admissible surface $\Xi$, a polygon $K$ and an asymptotically rigid homeomorphism $\psi$ such that $x=[\Xi,\psi]$ and $y= [\Xi \cup K, \psi]$. Since we also have $x=[\Sigma, \varphi]$, we know that $\varphi^{-1}\psi$ is isotopic to a homeomorphism $\xi$ that maps $\Xi$ to $\Sigma$ and that is rigid outside $\Xi$. Notice that $\xi$ maps $K$ to a polygon $H$ adjacent to $\Sigma$, so $\xi$ sends $\Xi \cup K$ to $\Sigma \cup H$ and is rigid outside $\Xi \cup K$, hence $y = [\Sigma \cup H, \varphi]$ as desired.
\end{proof}

\noindent
We begin by observing that our complex is connected:

\begin{lemma}\label{lem:Connected}
The cube complex $\mathscr{C}(A)$ is connected.
\end{lemma}

\begin{proof}
Fix two vertices $[\Sigma_1, \varphi_1]$ and $[\Sigma_2, \varphi_2]$. Since $\varphi_1$ and $\varphi_2$ are asymptotically rigid, so is $\varphi_2^{-1} \varphi_1$.  Let $\Xi$ be a support of $\varphi_2^{-1} \varphi_1$ and let $\Sigma_1^+$ be an admissible subsurface containing $\Sigma_1 \cup \varphi_1^{-1} \varphi_2(\Sigma_2) \cup \Xi.$ By construction, $\varphi_2^{-1} \varphi_1$ is rigid outside $\Sigma_1^+$.  Set $\Sigma_2^+= \varphi_2^{-1} \varphi_1( \Sigma_1^+)$. 

\medskip \noindent
Because $\Sigma_1 \subset \Sigma_1^+$ and $\Sigma_2 \subset \Sigma_2^+$, adding polygons to $\Sigma_1$ and $\Sigma_2$ one by one produces a path in $\mathscr{C}(A)$ from $[\Sigma_1,\varphi_1]$ to $[\Sigma_1^+,\varphi_1]$ and a path from $[\Sigma_2, \varphi_2]$ to $[\Sigma_2^+, \varphi_2]$. By construction, $\varphi_2^{-1} \varphi_1$ maps $\Sigma_1^+$ to $\Sigma_2^+$ and is rigid outside $\Xi$. As $\Xi \subset \Sigma_1^+$, $\varphi_2^{-1} \varphi_1$ must be rigid outside $\Sigma_1^+$. Consequently $[\Sigma_1^+, \varphi_1]=[\Sigma_2^+, \varphi_2]$ and this concludes the proof.
\end{proof}

\begin{proof}[Proof of Theorem \ref{thm:contractible}.]
We begin by constructing specific contractible compact subcomplexes in $\mathscr{C}(A)$. First of all, we need to introduce some terminology. Given a vertex $x$ and an admissible subsurface $\Sigma$, we say that $[\Sigma,\mathrm{id}]$ \emph{dominates} $x$ if there exists an \emph{increasing path} from $x$ to $[\Sigma, \mathrm{id}]$, i.e. a path along which the height of a vertex is always greater than the height of the previous one. Given a finite collection of vertices $\mathcal{S}$ and an admissible subsurface $\Sigma$ such that $[\Sigma, \mathrm{id}]$ dominates every vertex in $\mathcal{S}$, we denote by $X(\mathcal{S},\Sigma)$ the subcomplex spanned by the union of all the increasing paths from a vertex in $\mathcal{S}$ to $[\Sigma, \mathrm{id}]$.

\medskip \noindent
Observe that dominating vertices always exist:

\begin{claim}\label{claim:Dom}
Let $\mathcal{S}$ be a finite collection of vertices. There exists an admissible subsurface $\Sigma$ such that $[\Sigma, \mathrm{id}]$ dominates all the vertices in $\mathcal{S}$.
\end{claim}

\noindent
Fix an enumeration $\mathcal{S}= \{ [\Sigma_i, \varphi_i] \mid i \in I\}$. For every $i\in I$, let $\Lambda_i$ be a support of $\varphi_i$ and let $\Omega_i$ denote the smallest admissible subsurface containing $\Sigma_i\cup \Lambda_i$. Fix an admissible subsurface $\Sigma$ that contains $\varphi_i(\Omega_i)$ for every $i \in I$. Adding polygons to $\Sigma_i$ produces an increasing path from $[\Sigma_i, \varphi_i]$ to $[\Omega_i,\varphi_i]= [\varphi_i(\Omega_i), \mathrm{id}]$; and next adding polygons to $\varphi_i(\Omega_i)$ produces an increasing path from $[\varphi_i(\Omega_i),\mathrm{id}]$ to $[\Sigma, \mathrm{id}]$. Thus, we have proved that $[\Sigma, \mathrm{id}]$ dominates each vertex in $\mathcal{S}$, concluding the proof of our claim.

\medskip \noindent
Our goal now is to show that the subcomplexes $X(\mathcal{S}, \Sigma)$ are always contractible.

\medskip \noindent
So fix a finite collection of vertices $\mathcal{S}$ and an admissible subsurface $\Sigma$ such that $[\Sigma, \mathrm{id}]$ dominates all the vertices in $\mathcal{S}$. 
We assume that $X(\mathcal{S},\Sigma)$ is not reduced to a single vertex, i.e. $X(\mathcal{S}, \Sigma) \neq \{[\Sigma, \mathrm{id}] \}$. As a consequence, a vertex $x \in \mathcal{S}$ of minimal height admits neighbors in $X(\mathcal{S}, \Sigma)$, say $x_1, \ldots,x_k$ where $k\geq 1$. Notice that $x$ must also have minimal height in $X(\mathcal{S},\Sigma)$, so it follows from Lemma \ref{lem:Choice} that there exist an admissible subsurface $\Xi$, pairwise distinct polygons $H_1, \ldots, H_k$ and a homeomorphism $\varphi$ such that $x= [\Xi, \varphi]$ and $x_i= [\Xi \cup H_i, \varphi]$ for every $1 \leq i \leq k$. 

\begin{claim}\label{claim:CubeinX}
The cube spanned by $\left\{ \left[\Xi \cup \bigcup\limits_{i \in I} H_i, \varphi \right] \mid I \subset \{1, \ldots, k\} \text{ finite} \right\}$ lies in $X(\mathcal{S},\Sigma)$. 
\end{claim}

\noindent
Fix a non-empty finite subset $I \subset \{1, \ldots, k\}$. Because there exists an increasing path from $x$ to $[\Sigma, \mathrm{id}]$, we know from Lemma \ref{lem:Choice} that there exist polygons $P_1, \ldots, P_r$ such that $[\Xi \cup P_1 \cup \cdots \cup P_r, \varphi]= [\Sigma,\mathrm{id}]$. Given an $i \in I$, notice that $x_i$ lies on an increasing path from $x$ to $[\Sigma, \mathrm{id}]$, so we deduce again from Lemma \ref{lem:Choice} that there exist polygons $Q_1, \ldots, Q_s$ such that $[\Sigma, \mathrm{id}]=[\Xi \cup H_i \cup Q_1 \cup \dots \cup Q_s, \varphi]$. From
$$[\Xi \cup P_1 \cup \cdots \cup P_r, \varphi]= [\Sigma,\mathrm{id}]= [\Xi \cup H_i \cup Q_1 \cup \dots \cup Q_s, \varphi]$$
it follows that $\Xi \cup P_1 \cup \cdots \cup P_r = \Xi \cup H_i \cup Q_1 \cup \dots \cup Q_s$. Thus, we have proved that, for every $i \in I$, there exists $1 \leq \sigma(i) \leq r$ such that $H_i= P_{\sigma_i}$. Because the $H_i$'s are adjacent to $\Xi$, adding to $\Xi$ the $H_i$'s and next the remaining $P_i$'s produces an increasing path from $x$ to $[\Sigma, \mathrm{id}]$ passing through $\left[\Xi \cup \bigcup\limits_{i \in I} H_i, \varphi \right]$, concluding the proof of our claim.

\medskip \noindent
We deduce from Claim \ref{claim:CubeinX} that $X(\mathcal{S}, \Sigma)$ deformation retracts onto the proper subcomplex $X\left( (\mathcal{S}\backslash \{x\}) \cup  \{x_1, \ldots, x_r\}, \Sigma \right)$. Moreover, notice that $[\Sigma,\mathrm{id}]$ dominates all the vertices in $(\mathcal{S}\backslash \{x\}) \cup  \{x_1, \ldots, x_r\}$. Therefore, by iterating the process, we find a sequence of finite sets $\mathcal{S}_1, \mathcal{S}_2, \ldots$ of vertices all dominated by $[\Sigma,\mathrm{id}]$ such that
$$X(\mathcal{S},\Sigma) \supsetneq X(\mathcal{S}_1,\Sigma) \supsetneq X(\mathcal{S}_2, \Sigma) \supsetneq \cdots.$$
Since all these subcomplexes contain only finitely many cells, the sequence must eventually stop, and we find a finite collection $\mathcal{R}$ of vertices dominated by $[\Sigma,\id]$ such that $X(\mathcal{S}, \Sigma)$ deformation retracts onto $X(\mathcal{R}, \Sigma)$ and such that the previous process cannot apply to $X(\mathcal{R},\Sigma)$; in other words, $X(\mathcal{R}, \Sigma)= \{ [\Sigma, \mathrm{id}] \}$.
Thus, we have proved that $X(\mathcal{S}, \Sigma)$ is contractible. 

\medskip \noindent
We are finally ready to prove our proposition. Given an $n \geq 1$ and a continuous map $f : \mathbb{S}^n \to \mathscr{C}(A)$, the image of $f$ lies in a compact subcomplex, say $Y$. According to Claim~\ref{claim:Dom}, there exists an admissible subsurface $\Sigma$ such that $[\Sigma, \mathrm{id}]$ dominates each vertex in $Y$, so $Y$ lies in the contractible subcomplex $X(Y^{(0)}, \Sigma)$, proving that $f$ is homotopically trivial. Thus, all the homotopy groups of $\mathscr{C}(A)$ are trivial, which implies that our cube complex is contractible thanks to Whitehead's theorem.
\end{proof}

\subsection{A word about curvature}\label{section:curvature}

\noindent
In view of Theorem \ref{thm:contractible}, it is natural to ask whether our cube complex is nonpositively curved. Recall that a cube complex is \emph{CAT(0)} if it is simply connected and if the links of its vertices are simplicial flag complexes. 

\begin{figure}
\begin{center}
\includegraphics[scale=0.27]{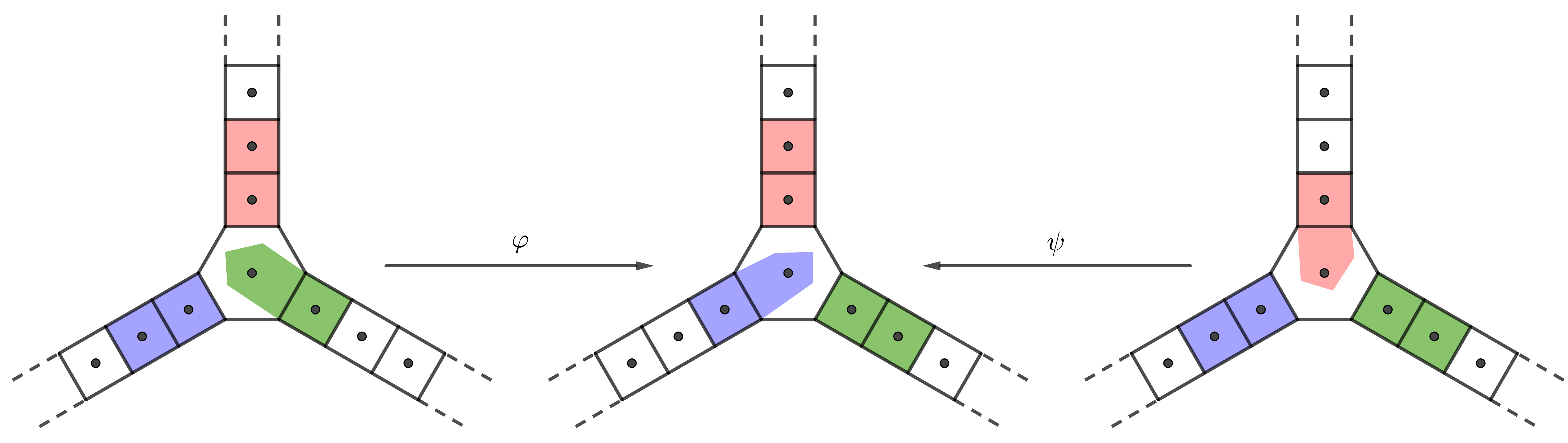}
\caption{The homeomorphisms $\varphi$ and $\psi$.}
\label{phipsi}
\end{center}
\end{figure}

\begin{figure}
\begin{center}
\includegraphics[scale=0.35]{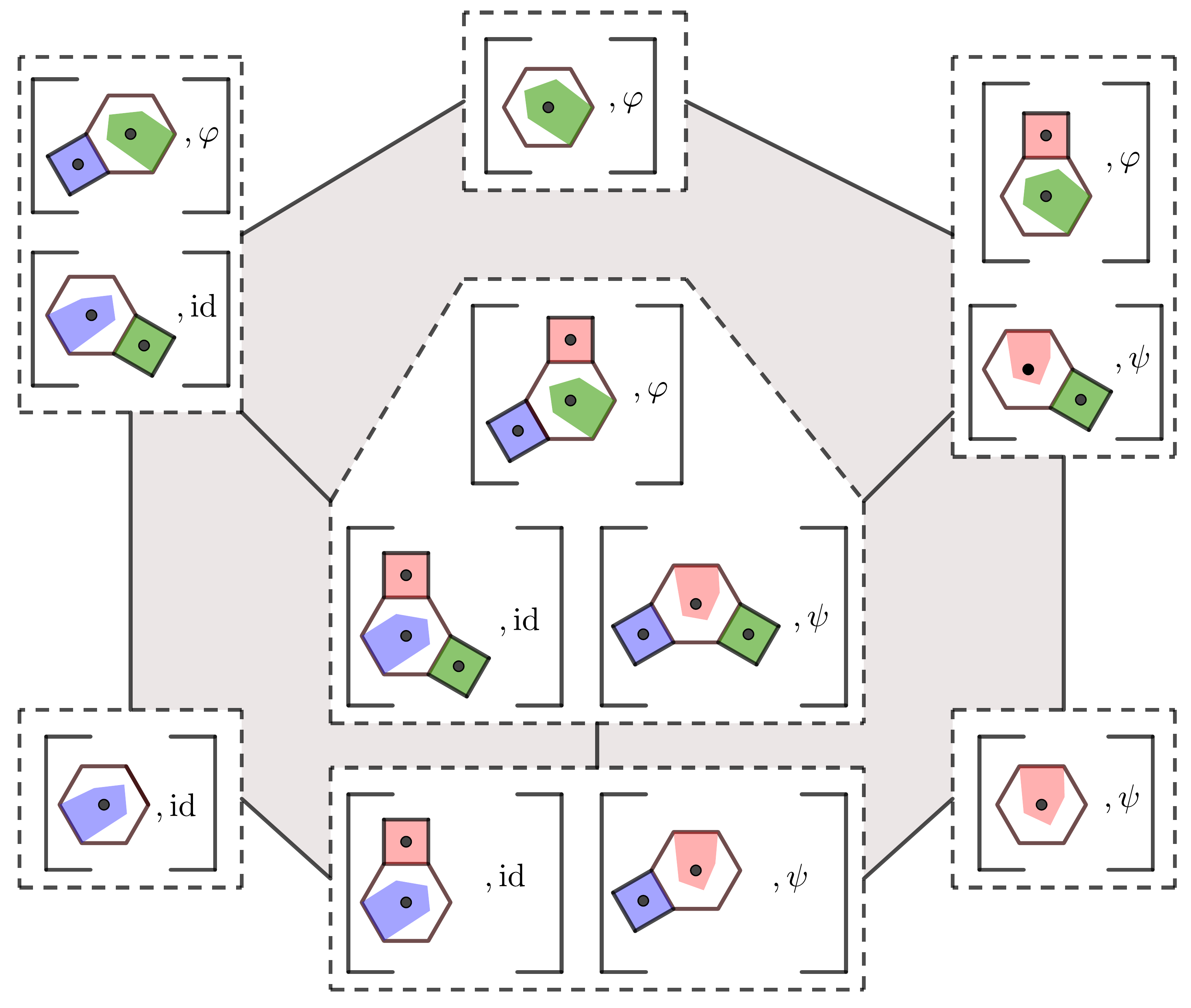}
\caption{}
\label{noCAT}
\end{center}
\end{figure}

\medskip \noindent
Unfortunately, our cube complexes may not be CAT(0). As an example, let $R_3$ denote the union of three infinite rays with a common origin and let $\varphi,\psi$ be the two homeomorphisms illustrated by Figure \ref{phipsi}. Then Figure \ref{noCAT} shows three squares in $\mathscr{C}(R_3)$ pairwise intersecting along an edge and all three intersecting along a vertex. However, this subcomplex does not span a $3$-cube as the missing vertex would have to be of height zero, which is impossible. Thus, the cube complex $\mathscr{C}(R_3)$ contains a vertex whose link is not flag, and a fortiori is not CAT(0).

\medskip \noindent
Nevertheless, our cube complexes turn out to be CAT(0) in several cases of interest, including the braided Ptolemy-Thompson groups. And, in other cases, the construction can be adaptated. This will be explained in a forthcoming article \cite{Next}.

\section{Classification of finite subgroups}\label{section:torsion}

\noindent
In this section, we fix a locally finite tree $A$ embedded into the plane in such a way that its vertex-set is discrete. As a first application of the cube complex constructed in the previous section, let us notice that we can classify the finite subgroups in $\amod(A)$. The key observation is the following:

\begin{lemma}\label{prop:finitegrps}
Let $G \leq \amod(A)$ be a finite subgroup. Then $G$ fixes a vertex of $\mathscr{C}(A)$ of the form $[\Sigma,\id]$.
\end{lemma}

\begin{proof}
Let $\Sigma' \subset \mathscr{S}^\sharp(A)$ be an admissible subsurface that is a support for all elements $g\in G$. Such a surface exists since $G$ is finite, by assumption. Let $\Sigma \subset \mathscr{S}^\sharp(A)$ be the smallest admissible subsurface containing the admissible subsurfaces $\{g(\Sigma')\mid g\in G\}$. We claim that the vertex $[\Sigma,\id]$ is fixed by $G$.

\medskip \noindent
So let $k \in G$ be an arbitrary element. Then $k(\Sigma)$ is again admissible and $k$ is rigid outside $\Sigma$, since $\Sigma$ is a support of $k$. Moreover, $\{g(\Sigma')\mid g\in G\}$ is invariant under $k$, so we must have $k(\Sigma)=\Sigma$. Consequently, we have
$$[\Sigma, \id]=[\Sigma, k]= k \cdot [\Sigma, \id],$$ 
i.e. $k$ fixes $[\Sigma, \id]$, as desired.
\end{proof}

\noindent
Therefore, the classification of finite subgroups in $\amod(A)$ reduces to the classification of finite subgroups in vertex-stabilisers of $\mathscr{C}(A)$. The structure of these stabilisers are described by the following statement:

\begin{lemma}\label{lem:VertexStab}
The stabiliser in $\amod(A)$ of a vertex $[\Sigma, \id]$ in $\mathscr{C}(A)$ is a subgroup of $\mathrm{stab}(\Sigma)$ in $\mathrm{Mod}(\mathscr{S}^\sharp(A))$, and it satisfies
$$1 \to \mathrm{Mod}(\Sigma) \to \stab([\Sigma, \id]) \to \mathbb{Z}_{r(\Sigma)} \to 1$$
for some integer $r(\Sigma) \geq 0$, where the morphism to $\mathbb{Z}_{r(\Sigma)}$ comes from the action by cyclic permutations of $\mathrm{stab}([\Sigma,\id])$ on components of $\mathrm{Fr}(\Sigma)$. 
\end{lemma}

\noindent
For instance, if $A$ is the $n$-regular tree and if $\Sigma$ is a single polygon, then $r(\Sigma)=n$. However, if $A=R_n$ with $n \geq 3$ and if $\Sigma$ is again a single polygon, then $r(\Sigma)=n$ if $\Sigma$ is the central polygon and $r(\Sigma)=0$ otherwise.

\begin{proof}[Proof of Lemma \ref{lem:VertexStab}.]
The equivalence relation in the definition of vertices in $\mathscr{C}(A)$ implies directly that the stabiliser $\stab([\Sigma, \id])$ is the subgroup 
$$\left\{ g \in \amod(A) \mid \text{up to isotopy, $g$ maps $\Sigma$ to itself and is rigid outside $\Sigma$} \right\}.$$
Now, note that $\Mod(\Sigma)$ injects into the stabiliser $\stab([\Sigma, \id])$ in the following way: an isotopy class of homeomorphisms in $\Mod(\Sigma)$ represented by $g$ defines an element in $\amod(A)$ as the isotopy class of the asymptotically rigid element whose restriction to $\Sigma$ is $g$ and that is the identity outside $\Sigma$. 
Consider the following short exact sequence: 
\[ 1\longrightarrow \Mod(\Sigma)\overset{\iota}{ \hooklongrightarrow} \stab([\Sigma, \id]) \overset{\pi}{\longrightarrow} \Z_{r(\Sigma)} \longrightarrow 1, \] where $\iota$ is the inclusion and $\pi$ is defined as follows. If $g\in \stab([\Sigma, \id])$ then $g$ (up to isotopy) is an orientation-preserving homeomorphism of $\Sigma$. Consequently, it preserves the cyclic order of the connected components of $\mathrm{Fr}(\Sigma)$. This yields a homomorphism from the stabiliser $\stab([\Sigma, \id])$ to some finite cyclic group, whose image is hence a cyclic group $\mathbb{Z}_{r(\Sigma)}$.
\end{proof}

\noindent
We are now ready to state the main result of this section:

\begin{thm}\label{thm:FiniteOrders}
Every finite subgroup in $\amod(A)$ is cyclic. Moreover, the positive divisors of the integers in
$$\{ \gcd(r(\Sigma),h(\Sigma)), \ \gcd(r(\Sigma),h(\Sigma)-1) \mid \text{$\Sigma$ admissible, $r(\Sigma) \neq 0$}\},$$ where $r(\Sigma)$ is defined in Lemma \ref{lem:VertexStab}, are exactly the orders of the finite-order elements in $\amod(A)$.
\end{thm}

\noindent
Before turning to the proof of Theorem \ref{thm:FiniteOrders}, recall that the center of the braid group $\mathrm{Mod}(D_k)$ on $k$ strands is the infinite cyclic group generated by the \emph{full twist} $\tau$, which is a Dehn twist around a simple closed curve lying in a small tubular neighborhood of the boundary $\partial D_k$ and parallel to $\partial D_k$. An element $g\in \Mod(D_k)$ is \emph{periodic} if $g^r=\tau^s$ for some integers $r,s \in \mathbb{Z}$, not both zero. For instance, if we fix a cyclic order on the punctures of $D_k$, then the (class of the) homeomorphism $\varepsilon$ given by permuting the $k$ punctures cyclically is periodic (see Figure~\ref{fig:braid1}); and, if we fix a cyclic order on all the punctures of $D_k$ but one, then the (class of the) homeomorphism $\delta$ given by permuting cyclically the ordered punctures is also periodic (see Figure~\ref{fig:braid2}). As observed in \cite[Corollary~2.2]{MR2253663}, a consequence of a theorem proved by Ker\'ekj\'art\'o, Brouwer, and Eilenberg \cite{Eilenberg}, which we record for future use, is that $\varepsilon$ and $\delta$ are essentially the only periodic elements in our braid group:

\begin{figure}
	\centering
	\begin{minipage}{0.45\textwidth}
	\begin{center}

	\begin{tikzpicture}
	\draw[thick, ->] (2,0) arc (0:50:2);
	\draw[thick, ->] (0.6,1.8) arc (72:122:2);
	\draw[thick, ->] (-1.6,1.218) arc (144:194:2);
	\draw[thick, ->] (-1.6,-1.18) arc (216:266:2);
	\draw[thick, ->] (0.6,-1.8) arc (288:338:2);
	\filldraw[black] (1.9,-0.3) circle (2pt);
	\filldraw[black] (0.88,1.7) circle (2pt);
	\filldraw[black] (-1.4,1.5) circle (2pt);
	\filldraw[black] (-1.8,-0.8) circle (2pt);
	\filldraw[black] (0.2,-1.9) circle (2pt);
	\end{tikzpicture}
\end{center}
\caption{The element $\varepsilon$ in $\Mod(D_k)$.}
\label{fig:braid1}
	\end{minipage}\hfill
	\begin{minipage}{0.45\textwidth}
		\begin{center}
			
			\begin{tikzpicture}
			\draw[thick, ->] (2,0) arc (0:75:2);
			\draw[thick, ->] (0,2) arc (90:165:2);
			\draw[thick, ->] (-2,0) arc (180:255:2);
			\draw[thick, ->] (0,-2) arc (270:345:2);
			\filldraw[black] (1.98,-0.2) circle (2pt);
			\filldraw[black] (0.16,1.98) circle (2pt);
			\filldraw[black] (-1.98,0.2) circle (2pt);
			\filldraw[black] (-0.16,-1.98) circle (2pt);
			\filldraw[black] (0,0) circle (2pt);
			\end{tikzpicture}
			
		\end{center}
		\caption{The element $\delta$ in $\Mod(D_k)$.}
		\label{fig:braid2}
	\end{minipage}
\end{figure}

\begin{prop}\label{prop:periodic}
Every periodic element in $\Mod(D_k)$ is conjugate to a power of $\varepsilon$ or~$\delta$.
\end{prop}

\noindent
Now we are ready to prove our main theorem.

\begin{proof}[Proof of Theorem \ref{thm:FiniteOrders}.]
Let $G \leq \amod(A)$ be a finite subgroup. By Proposition~\ref{prop:finitegrps}, $G$ fixes a vertex $[\Sigma, \id]$. Recall from Lemma~\ref{lem:VertexStab} that there exists a short exact sequence
	\begin{equation}\label{eq:stab}
1\longrightarrow \Mod(\Sigma) \hooklongrightarrow \stab([\Sigma, \id]) \overset{\pi}{\longrightarrow} \Z_{r(\Sigma)} \longrightarrow 1.
	\end{equation}
The subgroup $\mathrm{Mod}(\Sigma)$ is isomorphic to the mapping class group $\Mod(D_{h(\Sigma)})$ of the $h(\Sigma)$-punctured disk $D_{h(\Sigma)}$, itself isomorphic to the braid group on $h(\Sigma)$ strands. Since braid groups are torsion-free, the homomorphism $\pi$ injects $G$ into $\mathbb{Z}_{r(\Sigma)}$. In particular, $G$ is cyclic, proving the first assertion of our theorem. Moreover, if $r(\Sigma)=0$ then $G$ must be trivial, so, from now on, we assume that $r(\Sigma) \neq 0$.

\medskip \noindent
If $h(\Sigma)=1$, then $\mathrm{Mod}(\Sigma)=\{1\}$, which implies that $\mathrm{stab}([\Sigma,\id])$ is cyclic of order $r(\Sigma)$. From now on, we assume that $h(\Sigma) \geq 2$.

\medskip \noindent
Fix an element $\rho\in  \stab([\Sigma, \id])$ such that its restriction to some disk in $\Sigma$ containing all the punctures of $\Sigma$ is the identity and such that it cyclically shifts the components of $\mathrm{Fr}(\Sigma)$ with the smallest possible angle. (Here, following Lemma \ref{lem:VertexStab}, we are thinking of $\mathrm{stab}([\Sigma,\id])$ as a subgroup of $\mathrm{stab}(\Sigma) \leq \mathrm{Mod}(\mathscr{S}^\sharp(A))$, and so $\rho$ as a homeomorphism stabilising $\Sigma$.) Notice that $\rho$ is sent to a generator in $\mathbb{Z}_{r(\Sigma)}$; that $\rho^{r(\Sigma)}$ coincides with the full twist $\tau$ of $\Sigma$; that $\rho$ commutes with all the elements in $\mathrm{Mod}(\Sigma) \leq \mathrm{stab}([\Sigma,\id])$; and that $\rho$ has infinite order because $h(\Sigma) \geq 2$. 

\medskip \noindent
Let $g \in \mathrm{stab}([\Sigma,\id])$ be a non-trivial element. It can be written as $\rho^t \sigma$ where $t \in \mathbb{Z}$ and $\sigma \in \mathrm{Mod}(\Sigma)$; up to replacing $g$ with $g^{-1}$, we can suppose without loss of generality that $t \geq 0$. Notice that, if $\sigma$ is not periodic, then $g$ has infinite order. Indeed, if $g$ has finite order then $\sigma$ must have a power that is a power of $\rho$, so $\sigma$ must have a power that is a power of the full twist $\tau$. From now on, we assume that $\sigma$ is periodic. According to Proposition \ref{prop:periodic}, $g$ must be conjugate to $\rho^t \epsilon^{s}$ or $\rho^t \delta^{s}$ for some $s \in \mathbb{Z}$. Because $g$ is non-trivial, $t$ and $s$ cannot be both zero. Therefore, if $t = 0$ or $s=0$, then $g$ has infinite order. From now on, we assume that $t>0$ and $s \neq 0$.  

\medskip \noindent
First, assume that $g$ is conjugate to $\rho^t \epsilon^{s}$. Notice that $(t \vee r(\Sigma))/t$ is the smallest integer $n \geq 1$ such that $(\rho^t)^n$ belongs to $\langle \tau \rangle$, and similarly $(|s| \vee h(\Sigma))/|s|$ is the smallest integer $n \geq 1$ such that $(\epsilon^{s})^n$ belongs to $\langle \tau \rangle$. (Here, $\vee$ refers to the least common multiple.) Therefore, $\nu:= \frac{t \vee r(\Sigma)}{t} \vee \frac{|s| \vee h(\Sigma)}{|s|}$ is the smallest integer $n \geq 1$ such that $(\rho^t\epsilon^{s})^n$ belongs to $\langle \tau \rangle$. Notice that
$$\left( \rho^t\epsilon^{s} \right)^\nu = \left( \rho^{r(\Sigma)} \right)^{t \nu/ r(\Sigma)} \left( \epsilon^{h(\Sigma)} \right)^{s \nu/h(\Sigma)} = \tau^{\frac{t\nu}{r(\Sigma)} +\frac{s \nu}{h(\Sigma)}}.$$
Therefore, either $th(\Sigma) + s r(\Sigma) \neq 0$ and $g$ has infinite order (since it is conjugate to a non-trivial power of $\tau$) or $th(\Sigma) + s r(\Sigma)=0$ and $g$ has order $\nu$. 

\medskip \noindent
Similarly, if $g$ is conjugate to $\rho^t \delta^s$, we set $\nu:= \frac{t \vee r(\Sigma)}{t} \vee \frac{|s| \vee (h(\Sigma)-1)}{|s|}$ and we show that either $t(h(\Sigma)-1) + sr(\Sigma) \neq 0$ and $g$ has infinite order or $t(h(\Sigma)-1) + sr(\Sigma) = 0$ and $g$ has order $\nu$. 

\medskip \noindent
So far, we have proved that
$$\begin{array}{l} \displaystyle \left\{ \text{divisors of } r(\Sigma) \mid \Sigma \text{ single polygon, $r(\Sigma) \neq 0$} \right\} \\ \\ \displaystyle \cup \left\{  \frac{t \vee r(\Sigma)}{t} \vee \frac{s \vee h(\Sigma)}{s} \mid t,s>0, th(\Sigma) = s r(\Sigma), \text{$\Sigma$ admissible}, r(\Sigma) \neq 0 \right\} \\ \\ \displaystyle \cup \left\{  \frac{t \vee r(\Sigma)}{t} \vee \frac{s \vee (h(\Sigma)-1)}{s} \mid t,s>0, t(h(\Sigma)-1) = s r(\Sigma), \text{$\Sigma$ admissible}, r(\Sigma) \neq 0 \right\} \end{array}$$
is the set of orders of the finite-order elements in $\amod(A)$; according to our next claim, it coincides with the positive divisors of the integers in
$$\begin{array}{l} \displaystyle \left\{ r(\Sigma) \mid \text{$\Sigma$ single polygon, $r(\Sigma) \neq 0$} \right\} \\ \hspace{2cm} \displaystyle \cup \left\{ \gcd(r(\Sigma), h(\Sigma)), \gcd(r(\Sigma),h(\Sigma)-1) \mid \text{$\Sigma$ admissible, $r(\Sigma) \neq 0$} \right\}. \end{array}$$
We conclude the proof of our theorem by noticing that the left set lies in the right one.

\begin{claim}
For all integers $a,b \geq 1$, 
$$\left\{ \frac{t \vee a}{t} \vee \frac{s\vee b}{s} \mid t,s \geq 1, \ tb=sa \right\}$$
coincides with the set of positive divisors of $\mathrm{gcd}(a,b)$. 
\end{claim}

\noindent
First, assume that $t,s \geq 1$ are two integers satisfying $tb=sa$. Fix a prime $p$ and let 
$$v_p : n \mapsto \max\{k \geq 0 \mid p^k \text{ divides } n \}$$ 
denote the $p$-adic valuation. We distinguish two cases. If $v_p(t) \geq v_p(a)$, then
$$v_p \left( \frac{t \vee a}{t} \right) = \max(v_p(t),v_p(a))-v_p(t) =0.$$
Notice that $v_p(t)+v_p(b)=v_p(tb)=v_p(sa)=v_p(s)+v_p(a)$ implies that $v_p(s) \geq v_p(b)$, so we also have
$$v_p \left( \frac{s \vee b}{s} \right) = \max(v_p(s),v_p(b))-v_p(s)=0.$$
Therefore,
$$v_p \left( \frac{t \vee a}{t} \vee \frac{s \vee b}{s} \right) =\max \left( v_p \left( \frac{t \vee a}{t} \right), v_p \left( \frac{s \vee b}{s} \right) \right)=0.$$
Next, if $v_p(t)< v_p(a)$, then
$$v_p \left( \frac{t \vee a}{t} \right) = \max(v_p(t),v_p(a))-v_p(t) = v_p(a)-v_p(t).$$
Notice that $v_p(t)+v_p(b)=v_p(tb)=v_p(sa)=v_p(s)+v_p(a)$ implies that $v_p(s) < v_p(b)$, so we also have
$$v_p \left( \frac{s \vee b}{s} \right) = \max(v_p(s),v_p(b))-v_p(s)= v_p(b)-v_p(s).$$
Therefore,
$$v_p \left( \frac{t \vee a}{t} \vee \frac{s \vee b}{s} \right) =\max \left( v_p \left( \frac{t \vee a}{t} \right), v_p \left( \frac{s \vee b}{s} \right) \right)$$
coincides with $v_p(a)-v_p(t)= v_p(b)-v_p(s)$. We conclude that
$$v_p \left( \frac{t \vee a}{t} \vee \frac{s \vee b}{s} \right) \leq \min(v_p(a),v_p(b)) \text{ for every prime $p$},$$
proving that $\frac{t \vee a}{t} \vee \frac{s \vee b}{s}$ divides $\gcd(a,b)$. Thus, the integers in the set given by our claim are all positive divisors of $\gcd(a,b)$. Conversely, let $k$ be a positive divisor of $\gcd(a,b)$. Set $t:=a/k$ and $s:=b/k$, and notice that $tb= ab/k = sa$ and that
$$ \frac{t \vee a}{t} \vee \frac{s \vee b}{s}= \frac{a}{t} \vee \frac{b}{s} = k \vee k = k.$$
So the divisor $k$ belongs to the set given by our claim. This concludes the proof.
\end{proof}

\noindent
Let us illustrate Theorem \ref{thm:FiniteOrders} on two examples. First, we consider the trees associated to the braided Houghton groups.

\begin{prop}
Let $n \geq 2$ be an integer. Then $\amod(R_n)$ contains an element of order $l$ if and only if $l$ divides $n$.
\end{prop}

\begin{proof}
If $\Sigma$ is an admissible subsurface that does not contain the central vertex of $R_n$, then $r(\Sigma)=2$ if $n=2$ and $r(\Sigma)=0$ if $n \geq 3$ (because the complement of $\Sigma$ has two components, one of them one contains a $2n$-gon and the other contain only $4$-gons). If $\Sigma$ is an admissible subsurface containing the central vertex of $R_n$, then $r(\Sigma)=n$. The desired conclusion follows from Theorem \ref{thm:FiniteOrders}.
\end{proof}

\noindent
As a second example, we apply Theorem \ref{thm:FiniteOrders} to the braided Higman-Thompson groups. 

\torsionthm

\begin{proof}
First, assume that $m =n+1$. So $A_{n,m}$ is $(n+1)$-regular. If $\Sigma$ is an admissible subsurface, then 
$$r(\Sigma)= n+1 + (h(\Sigma)-1)(n-1) = 2 + h(\Sigma)(n-1).$$
Notice that $\gcd(r(\Sigma),h(\Sigma))$ must divide $r(\Sigma)-h(\Sigma)(n-1) = 2$ and that $\gcd(r(\Sigma),h(\Sigma)-1)$ must divide $r(\Sigma)-(h(\Sigma)-1)(n-1)=n+1$. Also, notice that, if $\Sigma$ has height $2$, then $\gcd(r(\Sigma),h(\Sigma))= \gcd(2n,2)=2$; and that, if $\Sigma$ has height $n+2$, then $\gcd(r(\Sigma),h(\Sigma)-1) = \gcd(n(n+1),n+1)=n+1$. Therefore, it follows from Theorem \ref{thm:FiniteOrders} that the possible finite orders in $\mathrm{br}T_{n,m}$ are the positive divisors of $2$ and $n+1$.

\medskip \noindent
From now on, we assume that $m \neq n+1$. If $\Sigma$ is an admissible subsurface that does not contain the vertex of valence $m$, then $r(\Sigma)=0$ because the complement of $\Sigma$ contains one component with a $2m$-gon while all the other components contain only $2(n+1)$-gons. But, if $\Sigma$ contains the vertex of valence $m$, then
$$r(\Sigma)= m + (h(\Sigma)-1)(n-1) = m-n+1 + h(\Sigma)(n-1).$$
Notice that $\gcd(r(\Sigma),h(\Sigma)-1)$ must divide $r(\Sigma)-(h(\Sigma)-1)(n-1)=m$ and that, if $\Sigma$ has height $m+1$, then $\gcd(r(\Sigma),h(\Sigma)-1))=\gcd(mn,m)=m$. We need to distinguish two cases:
\begin{itemize}
	\item Assume that $m \neq n-1$. Notice that $\gcd(r(\Sigma),h(\Sigma))$ must divide $r(\Sigma)-h(\Sigma)(n-1)=m-n+1$. Moreover, if $\Sigma$ has height $\lvert m-n+1\rvert$, then $r(\Sigma)=\lvert m-n+1\rvert(n-1+ \epsilon)$ where $\epsilon \in \{\pm 1\}$ according to the sign of $m-n+1$ and we have $\gcd(r(\Sigma), h(\Sigma))=|m-n+1|$. It follows from Theorem~\ref{thm:FiniteOrders} that the finite orders in $\mathrm{br}T_{n,m}$ are the positive divisors of the integers in $\{m\} \cup \{m-n+1\}.$
	\item Assume that $m=n-1$. Then $\gcd(r(\Sigma),h(\Sigma))= \gcd(h(\Sigma)(n-1),h(\Sigma))=h(\Sigma)$. It follows from Theorem~\ref{thm:FiniteOrders} that all the possible orders occur in $\mathrm{br}T_{n,m}$.
\end{itemize}
This concludes the proof of our theorem.\qedhere
\end{proof}

\noindent
Interestingly, Theorem \ref{prop:torsion} allows us to distinguish many braided Higman-Thompson groups up to isomorphism. For instance:

\begin{cor}
Let $m,n \geq 2$ be two integers. The braided Ptolemy-Thompson groups $\mathrm{br}T_n$ and $\mathrm{br}T_m$ are  isomorphic if and only if $m=n$. 
\end{cor}

\begin{remark}\label{rem:funar_kap_non-iso}
From Theorem~\ref{prop:torsion}, we can directly recover the result of Funar-Kapoudjian that $\amod^*(A_2) \simeq \mathrm{br}T_{2,4}$ is not isomorphic to $\amod(A_2)=\mathrm{br}T_{2,3}$. More generally, we see that $\amod^*(A_n)\simeq \mathrm{br}T_{n,2n}$ is not isomorphic to $\amod(A_{n})=\mathrm{br}T_{n,n+1}$.

\medskip \noindent
However, Theorem~\ref{prop:torsion} does not give us any information on whether $\mathrm{br}T_{6n-1}$ is isomorphic to $\mathrm{br}T_{6n-2,6n}$, for instance. Note also that an interesting phenomenon occurs for the braided Higman-Thompson groups of the form $\mathrm{br}T_{n,n-1}$ for $n\geq 2$. They contain an element of each order, and so they cannot be distinguished up to isomorphism either by using torsion.
\end{remark}

\section{Finiteness properties}\label{section:TypeF}

This section is dedicated to the proofs of Theorems \ref{PtolemyConnected} and~\ref{thm:brHfiniteness} from the introduction. Namely, we want the prove that the braided Higman-Thompson group $\mathrm{br}T_{n,m}$ is of type $F_\infty$ and that the braided Houghton group $\mathrm{br}H_n$ is of type $F_{n-1}$ but not of type $F_n$.

\medskip \noindent
If we allow $n=1$ in our notation $A_{n,m}$, then $\amod(A_{n,m})$ coincides with the braided Higman-Thompson group $\mathrm{br}T_{n,m}$ if $n \geq 2$, and it contains the braided Houghton group $\mathrm{br}H_n$ as a finite-index subgroup if $n =1$. This notation allows us to work with these two families of groups simultaneously in Sections \ref{section:spine} and \ref{section:Links}.

\medskip \noindent
Our strategy is the following. First, we extract from the cube complex $\mathscr{C}(A_{n,m})$, on which $\amod(A_{n,m})$ acts, a smaller cube complex $\mathscr{SC}(A_{n,m})$, referred to as the \emph{spine}. We show in Section \ref{section:spine} that $\mathscr{SC}(A_{n,m})$ is a $\amod(A_{n,m})$-invariant subcomplex, on which $\mathscr{C}(A_{n,m})$ deformation retracts. In particular, $\mathscr{SC}(A_{n,m})$ is also contractible and its vertex-stabilisers are also finite extensions of braid groups. Next, in Section \ref{section:Links}, we describe the descending links in $\mathscr{SC}(A_{n,m})$, with respect to the height function, as complexes of arcs on the disc. Finally, in Sections \ref{section:HigmanThompson} and \ref{section:Houghton}, we study the connectedness of these complexes and we deduce our main theorems from standard arguments of Morse theory.

\subsection{The spine of the cube complex}\label{section:spine}

\noindent
Typically, there are two different types of vertices in the cube complex $\mathscr{C}(A_{n,m})$: those associated to the admissible subsurfaces containing the \emph{central polygon} (i.e. the polygon containing the vertex of degree $m$), and those associated to the admissible subsurfaces that do not contain the central polygon. The descending links (with respect to the height function) of the latter vertices are non-simply connected graphs if $n=1$, which is an obstruction to prove finiteness properties.  This observation motivates the following definition:

\begin{definition}
Fix two integers $n,m \geq 1$. The \emph{spine} $\mathscr{SC}(A_{n,m})$ is the subcomplex of $\mathscr{C}(A_{n,m})$ generated by the vertices
$$\{ [\Sigma, \varphi] \mid \text{$\Sigma$ contains the central polygon} \},$$
where the central polygon refers to an arbitrary polygon we fix once for all if $m=n+1$ and otherwise to the polygon of $\mathscr{S}^\sharp(A_{n,m})$ that contains the puncture corresponding to the unique vertex of $A_{n,m}$ having valence $m$.
\end{definition}

\noindent
It is clear that $\amod(A_{n,m})$ stabilises the spine $\mathscr{SC}(A_{n,m})$. By reproducing the proof of Theorem~\ref{thm:contractible} word for word, it can be proved that the spine is also contractible. In fact, we can prove more:

\begin{prop}\label{prop:HomotopyEqui}
For all $n,m \geq 1$, the spine $\mathscr{SC}(A_{n,m}) \subset \mathscr{C}(A_{n,m})$ is a $\amod(A_{n,m})$-invariant subcomplex, on which $\mathscr{C}(A_{n,m})$ deformation retracts. In particular, $\mathscr{C}(A_{n,m})$ and $\mathscr{SC}(A_{n,m})$ are homotopy equivalent.
\end{prop}

\begin{proof}
If $m=n+1$, then $A_{n,m}$ is the $m$-regular tree and $\mathscr{SC}(A_{n,m})=\mathscr{C}(A_{n,m})$. Indeed, if $[\Sigma, \varphi]$ is a vertex of $\mathscr{C}(A_{n,n+1})$, then there exists a rigid homeomorphism $t$ associated to an isometry of $A_{n,m}$ such that $t \Sigma$ contains the central polygon. Then $[\Sigma,\varphi]=[t \Sigma, \varphi \circ t^{-1}] \in \mathscr{SC}(A_{n,n+1})$. From now on, we assume that $m \neq n+1$. 

\medskip \noindent
As a consequence of Whitehead's theorem, it is sufficient to show that the inclusion $\mathscr{SC}(A_{n,m}) \hookrightarrow \mathscr{C}(A_{n,m})$ induces an isomorphism on each homotopy group in order to deduce that $\mathscr{C}(A_{n,m})$ deformation retracts on $\mathscr{SC}(A_{n,m})$. This is a consequence of the following statement: every finite subcomplex $F \subset \mathscr{C}(A_{n,m})$ is contained in another subcomplex $F^+$ that deformation retracts on $F^+ \cap \mathscr{SC}(A_{n,m})$. We need some preliminary observations before proving this statement.

\medskip \noindent
Let $\Sigma$ be an admissible subsurface. We refer to the \emph{central height} of $\Sigma$ as the minimal number of polygons to add to $\Sigma$ in order to get an admissible subsurface containing the central polygon, and we denote by $\tau(\Sigma)$ the subsurface obtained from $\Sigma$ by adding the polygon adjacent to $\Sigma$ that is the closest to the central polygon. If $\Sigma$ already contains the central polygon, then it has zero central height and $\tau(\Sigma)= \Sigma$.

\begin{claim}\label{claim:CentralHeight}
If $(\Sigma_1,\varphi_1)$ and $(\Sigma_2,\varphi_2)$ are two representatives of a vertex in $\mathscr{C}(A_{n,m})$, then $\Sigma_1$ and $\Sigma_2$ have the same central height and $[\tau(\Sigma_1),\varphi_1]=[\tau(\Sigma_2),\varphi_2]$. 
\end{claim}

\noindent
Because having central height zero can be read from the height and the number of components in the frontier of the admissible subsurface under consideration, our claim is clear if the vertex belongs to the spine. From now on, we assume that it does not belong to the spine, i.e. $\Sigma_1$ and $\Sigma_2$ do not contain the central polygon. We know that $\varphi_2^{-1}\varphi_1$ is isotopic to a homeomorphism $\varphi$ that sends $\Sigma_1$ to $\Sigma_2$ and that is rigid outside $\Sigma_1$. Necessarily, $\varphi$ stabilises the central polygon. Also, it sends the polygons between $\Sigma_1$ and the central polygon to the polygons between $\Sigma_2$ and the central polygon; in particular, it sends $\tau(\Sigma_1)$ to $\tau(\Sigma_2)$. The latter assertion shows that $[\tau(\Sigma_1),\varphi_1]=[\tau(\Sigma_2),\varphi_2]$; and the former assertion shows that $\Sigma_1$ and $\Sigma_2$ have the same central height. Thus, our claim is proved.

\medskip \noindent
Claim \ref{claim:CentralHeight} allows us to define the central height of a vertex $[\Sigma,\varphi]$ as the central height of $\Sigma$, and to define the map 
$$\tau : \left\{ \begin{array}{ccc} \mathscr{C}(A_{n,m})^{(0)} & \to & \mathscr{C}(A_{n,m})^{(0)} \\ \left[ \Sigma, \varphi \right] & \mapsto & \left[ \tau(\Sigma), \varphi \right] \end{array} \right..$$
Notice that, for every vertex $x \in \mathscr{C}(A_{n,m})$, $\tau^k(x)$ belongs to $\mathscr{SC}(A_{n,m})$ if $k$ is at least the central height of $x$. Consequently, if $F$ is an arbitrary finite subcomplex in $\mathscr{C}(A_{n,m})$, then the subcomplex $F^+$ spanned by
$$\left\{ \tau^k(x) \mid x \in F^{(0)}, k \geq 0 \right\}$$
is also finite and it intersects $\mathscr{SC}(A_{n,m})$. Let $S \subset F^+$ denote the subcomplex generated by the vertices of maximal central height. Notice that, as a consequence of Claims~\ref{claim:Tau} and~\ref{claim:InjectiveComponent}, for each connected component $R$ of $S$ the subcomplex in $F^+$ generated by $R \cup \tau(R^{(0)})$ is isomorphic to $R \times [0,1]$, where $R$ is sent to $R \times \{0\}$ and $\tau(R^{(0)})$ to $R \times \{1\}$. Therefore, we can retract simultaneously all the components of $S$ onto the subcomplex of $F^+$ spanned by the vertices that do not have maximal central height. By iterating the process, we eventually deformation retract $F^+$ onto its subcomplex spanned by the vertices of central height zero, i.e. onto $F^+ \cap \mathscr{SC}(A_{n,m})$ as desired.

\begin{claim}\label{claim:Tau}
For a $k$-cube $C$ in $S$, $C \cup \tau(C^{(0)})$ spans a $(k+1)$-cube.
\end{claim}

\noindent
There exist an admissible subsurface $\Sigma$, a homeomorphism $\varphi$ and polygons $H_1, \ldots, H_k$ such that the vertices of $C$ are
$$\left\{ \left[ \Sigma\cup \bigcup\limits_{i \in I} H_i ,\varphi \right] \mid I \subset \{1 ,\ldots, k\} \right\}.$$
Because they all have maximal central height, the polygons $H_1, \ldots, H_k$ are distinct from the polygon adjacent to $\Sigma$ that is the closest to the central polygon. Let $H_{k+1}$ denote the former polygon. Then the vertices of $C \cup \tau(C^{(0)})$ are
$$\left\{ \left[ \Sigma\cup \bigcup\limits_{i \in I} H_i ,\varphi \right] \mid I \subset \{1 ,\ldots, k+1\} \right\}.$$ 
Thus, $C \cup \tau(C^{(0)})$ spans a $(k+1)$-cube.

\begin{claim}\label{claim:InjectiveComponent}
The map $\tau$ is injective on each connected component of $S$.
\end{claim}

\noindent
Before turning to the proof of our claim, we need to introduce some notation. Given an admissible subsurface $\Sigma$ that does not contain the central polygon, we denote by $\alpha(\Sigma)$ the polygon of $\Sigma$ that is the closest to the central polygon, and $\omega(\Sigma)$ the adjacent polygon of the central polygon that is the closest to $\Sigma$. 

\medskip \noindent
Now, fix two vertices $a,b \in S$ belonging to the same connected component and assume that $\tau(a)= \tau(b)$. Given a path of vertices $x_1, \ldots, x_k \in S$ from $a$ to $b$, we claim that they admit representatives $(\Sigma_1,\varphi_1), \ldots, (\Sigma_k,\varphi_k)$ such that, for every $1 \leq i \leq k-1$, we have $\alpha(\Sigma_i) = \alpha(\Sigma_{i+1})$  and $\varphi_{i+1}^{-1}\varphi_i$ is isotopic to an asymptotically rigid homeomorphism that is rigid on the connected subsurface delimited by $\alpha(\Sigma_i)$ and that contains the central polygon. We construct our representatives inductively, starting with an arbitrary representative $(\Sigma_1,\varphi_1)$ of $x_1$.

\medskip \noindent
Assume that $(\Sigma_i, \varphi_i)$ is defined for some $1 \leq i \leq k-1$ and that the height of $x_{i+1}$ is larger than the height of $x_i$. According to Lemma \ref{lem:Choice}, there exists a polygon $H$ such that $(\Sigma_{i}\cup H,\varphi_i)$ is a representative of $x_{i+1}$. Then set $\Sigma_{i+1}=\Sigma_i \cup H$ and $\varphi_{i+1}= \varphi_i$. Observe that $H$ does not separate $\alpha(\Sigma_i)$ and the central polygon because $x_i$ and $x_{i+1}$ have the same central height, so $\alpha(\Sigma_{i+1})= \alpha(\Sigma_i)$.

\medskip \noindent
Now, assume that $(\Sigma_i,\varphi_i)$ is defined for some $1 \leq i \leq k-1$ and that the height of $x_{i+1}$ is smaller than the height of $x_i$. Fix a representative $(\Sigma, \varphi)$ of $x_{i+1}$. Up to rotating around the central polygon by a rigid homeomorphism, we assume that $\omega(\Sigma)= \omega(\Sigma_i)$. If $\alpha(\Sigma_i) \neq \alpha(\Sigma)$, then we denote by $K$ the polygon separating both $\Sigma$ and $\Sigma_i$ from the central polygon and that is the farthest from the central polygon; see Figure \ref{Tau}. We also denote $P$ (resp. $Q$) the adjacent polygon of $K$ that is the closest to $\Sigma_i$ (resp. $\Sigma$).  Now, according to Lemma \ref{lem:Choice}, there exists a polygon $H$ such that $x_i=[\Sigma \cup H, \varphi]$. From the equality $[\Sigma \cup H, \varphi]= [\Sigma_i, \varphi_i]$, we know that $\varphi^{-1} \varphi_i$ is isotopic to a homeomorphism $\psi$ that sends $\Sigma_i$ to $\Sigma \cup H$ and that is rigid outside $\Sigma_i$. Necessarily, $\psi$ fixes $K$ and the central polygon, which are distinct. Consequently, $\psi$ must stabilise each component of the frontier of $K$, contradicting the fact that $\psi$ must also send $P$ to $Q$. Thus we have proved that $\alpha(\Sigma_i)= \alpha(\Sigma)$. Observe that, because $x_i$ and $x_{i+1}$ have the same central height, necessarily $\alpha(\Sigma)= \alpha(\Sigma \cup H)$. Therefore, we can set $\Sigma_{i+1}= \Sigma$ and $\varphi_{i+1}=\varphi$. 
\begin{figure}
\begin{center}
\includegraphics[scale=0.3]{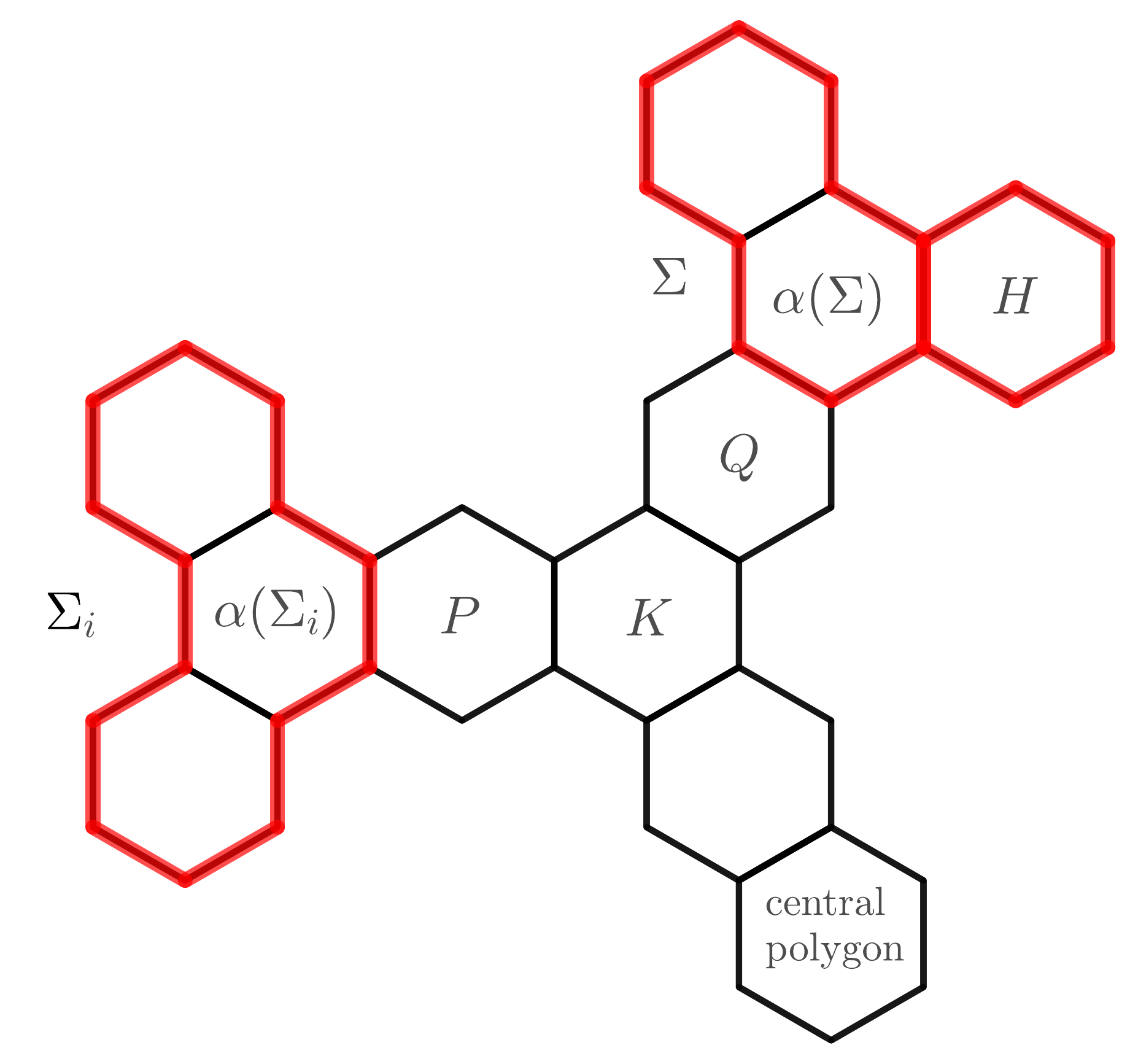}
\caption{}
\label{Tau}
\end{center}
\end{figure}

\medskip \noindent
Thus, our representatives are constructed. Since $\alpha(\Sigma_i)= \alpha(\Sigma_{i+1})$ for every $1 \leq i \leq k-1$, we denote by $A$ this common polygon. Also, we denote by $B$ the adjacent polygon of $A$ that is the closest to the central polygon. By construction,
$$\varphi_k^{-1}\varphi_1 = (\varphi_k^{-1} \varphi_{k-1}) \cdot (\varphi_{k-1}^{-1} \varphi_{k-2}) \cdots (\varphi^{-1}_2 \varphi_1)$$
is isotopic to an asymptotically rigid homeomorphism that is rigid on the connected subsurface delimited by $A$ and containing the central polygon. Observe that such a homeomorphism must stabilise $B$. On the other hand, it follows from the equality
$$[\Sigma_1 \cup B , \varphi_1] = \tau(x_1)=\tau(a)=\tau(b)= \tau(x_k) = [\Sigma_k \cup B, \varphi_k]$$
that $\varphi_k^{-1}\varphi_1$ is isotopic to an asymptotically rigid homeomorphism sending $\Sigma_1 \cup B$ to $\Sigma_k \cup B$ that is rigid outside $\Sigma_1 \cup B$. We deduce that $\varphi_{k}^{-1} \varphi_1$ is isotopic to an asymptotically rigid homeomorphism that $\Sigma_1$ to $\Sigma_k$ and that is rigid outside $\Sigma_1$, hence $a=[\Sigma_1,\varphi_1]= [\Sigma_k, \varphi_k]=b$ concluding the proof of Claim \ref{claim:InjectiveComponent}.
\end{proof}

\subsection{Description of the descending links}\label{section:Links}

\noindent
Let us recall some basic definitions about links and Morse functions in cube complexes.

\begin{definition}
Let $X$ be a cube complex. A map $f : X \to \mathbb{R}$ is a \emph{Morse function} if it is affine and non-constant on each cube of positive dimension and if the image $f(X^{(0)})$ is discrete. For every $m \in \mathbb{R}$, the \emph{sublevel set} $X_m$ is the subcomplex of $X$ generated by the vertices in $y \in X$ satisfying $f(y) \leq m$. The \emph{link} of a vertex $x \in X$ is the simplicial complex whose vertices are the edges containing $x$ and whose simplices are collections of edges spanning cubes. The \emph{descending link} of $x \in X$ (with respect to $f$) is the link of $x$ in the $X_{f(x)}$.
\end{definition}

\noindent
For $\mathscr{SC}(A_{n,m})$, the affine extension of the height function is our Morse function. Our goal in this section is to describe the descending links with respect to this function.

\medskip \noindent
Let $p, q\geq1$ and $r \geq 0$ be three integers. Fix a disc $\mathbb{D}$ with $p$ punctures in its interior and $q$ marked points on its boundary. Let $\{m_i \mid i \in \mathbb{Z}_q \}$ denote these marked points, ordered cyclically. From now on, an arc in $\mathbb{D}$ will refer to an arc that starts from a marked point and that ends at a puncture. Two arcs starting from the marked points $m_i,m_j$ are \emph{$r$-separated} if they are disjoint and if the distance between $i$ and $j$ in $\mathbb{Z}_q$ is $>r$ (where $\mathbb{Z}_q$ is metrically thought of as the cycle $\mathrm{Cayl}(\mathbb{Z}_q,\{1\})$). Notice that being $0$-separated amounts to being disjoint. We define $\mathfrak{C}(p,q,r)$ as the simplicial complex whose vertices are the isotopy classes of arcs and whose simplices are collections of arcs that are pairwise $r$-separated (up to isotopy). Observe that pairwise $r$-separated classes of arcs $\alpha_1,\dots, \alpha_n$ can be represented by pairwise $r$-separated arcs $\beta_1,\dots, \beta_n$.

\medskip \noindent
The braid group $\mathrm{Mod}(\mathbb{D})$ naturally acts on $\mathfrak{C}(p,q,r)$. Moreover:

\begin{lemma}\label{lem:LinkTopo}
For every $i \in \mathbb{Z}_q$, fix an arc $\alpha_i$ from $m_i$ to a puncture and assume that, if two marked points $m_i,m_j$ satisfy $\left|  i-j\right|  >r$, then $\alpha_i$ and $\alpha_j$ are disjoint. The subcomplex spanned by $\{\alpha_i, i \in \mathbb{Z}_q\}$ is a strict fundamental domain for the action $\mathrm{Mod}(\mathbb{D}) \curvearrowright \mathfrak{C}(p,q,r)$. Moreover, for every $I \subset \mathbb{Z}_q$ such that $i$ and $j$ are at distance $>r$ for all distinct $i,j \in I$, the stabiliser of the simplex spanned by $\{\alpha_i, i \in I\}$ coincides with the subgroup $\mathrm{Mod} \left( \mathbb{D} \backslash \bigcup\limits_{i \in I} N(\alpha_i) \right)$ where $N(\alpha_i)$ is a small tubular neighborhood of $\alpha_i$ for every $i \in I$. 
\end{lemma}

\begin{proof}
Let $\mathscr{S} \subset \mathfrak{C}(p,q,r)$ denote the subcomplex spanned by $\{\alpha_0, \ldots, \alpha_{q-1}\}$. A simplex $\{\alpha_i, i \in I\}$ in $\mathscr{S}$, $I \subset \mathbb{Z}_q$, is uniquely determined by the marked points $\{m_i, i \in I\} \subset \partial \mathbb{D}$. Because $\mathrm{Mod}(\mathbb{D})$ fixes the boundary $\partial \mathbb{D}$ pointwise, it follows that no two subsimplices in $\mathscr{S}$ belong to the same $\mathrm{Mod}(\mathbb{D})$-orbit. Next, let $\mathscr{S}'$ be an arbitrary simplex in $\mathfrak{C}(p,q,r)$. The vertices of $\mathscr{S}'$ can be represented by pairwise $r$-separated arcs $\beta_1, \ldots, \beta_s$. For every $1 \leq i \leq s$, let $m_{\sigma(i)}$ denote the unique marked point that belongs to $\beta_i$. Then there exists a homeomorphism $g \in \mathrm{Mod}(\mathbb{D})$ such that $g \cdot \beta_i = \alpha_{\sigma(i)}$ for every $1 \leq i \leq s$. By construction, $g \mathscr{S}' \subset \mathscr{S}$. Thus, we have proved that $\mathscr{S}$ is a strict fundamental domain for the action $\mathrm{Mod}(\mathbb{D}) \curvearrowright \mathfrak{C}(p,q,r)$. The second assertion of our lemma is clear.
\end{proof}

\noindent
We are now ready to prove the main result of this section.

\begin{prop}\label{prop:deslink}
Let $n ,m \geq 1$ be two integers. In $\mathscr{SC}(A_{n,m})$, the descending link of a vertex of height $k \geq m+1$ is isomorphic to $\mathfrak{C}(k, m+(k-1)(n-1), n-1)$. 
\end{prop}

\begin{proof}
Fix a vertex $[\Sigma, \varphi] \in \mathscr{SC}(A_{n,m})$ of height $k \ge m+1$. Up to translating our vertex by $\varphi^{-1}$, we can assume without loss of generality that $\varphi=\mathrm{id}$. We denote by $\mathscr{A}$ the set of collections of $n$ consecutive frontier-arcs of $\Sigma$. For every collection of $n$ consecutive frontier-arcs $i \in \mathscr{A}$, we choose a topological disc $D_i \subset \Sigma$ and a puncture $p_i \in \Sigma$ such that {$D_i\cap \partial\Sigma$ is connected, such that}  $D_i$ intersects the frontier of $\Sigma$ exactly in $i$, and such that $p_i$ is the only puncture in $D_i$. Moreover, we assume that $D_i \cap D_j = \emptyset$ for all $i, j \in \mathscr{A}$ satisfying $i \cap j= \emptyset$, which is possible because $k \geq m+1$.

\medskip \noindent
For every set $I \subset \mathscr{A}$ of pairwise disjoint collections of frontier-arcs, we fix an admissible subsurface $\Sigma_I$ and an asymptotically rigid homeomorphism $\varphi_I$ of $\mathscr{S}^\sharp(A_{n,m})$ such that $\varphi_I(\Sigma_I)=\Sigma$, $\varphi_I$ is rigid outside $\Sigma_I$, and $\varphi_I^{-1}(D_i)$ is an extremal polygon $H_i^I \subset \Sigma_I$ of the rigid structure for every $i \in I$. Observe that these conditions imply that $[\Sigma,\mathrm{id}]= [\Sigma_I, \varphi_I]$. For instance, $\Sigma_I$ can be constructed as follows. Since $\Sigma$ has height $k$, its frontier contains $m+(k-1)(n-1)$ arcs, which implies that $n \#I \leq m+(k-1)(n-1)$. Because $k \geq m+1$, we deduce that $\# I \leq k-1$. Consequently, we can fix an admissible subsurface $\Xi$ of height $k-\#I \geq 1$ containing the central polygon. Notice that the frontier of $\Xi$ contains $\left[ m+(k-1)(n-1) \right] - (n-1) \#I\geq \# I$ arcs in its frontier, so we can fix a collection $J$ of $\#I$ arcs in the frontier of $\Xi$ and a bijection $\iota : J \to I$ so that any two consecutive arcs in $J$ (with respect to the cyclic order) are separated by the same number of arcs that separate their images under $\iota$ in the frontier of $\Sigma$. Then setting $\Sigma_I$ as the subsurface obtained from $\Xi$ by adding a polygon for each arc in $J$ provides the desired subsurface.

\medskip \noindent
If $I \subset \mathscr{A}$ is a set of pairwise disjoint collections,  we denote by $Q_I$ the cube generated by 
$$\left[ \Sigma_I \backslash \bigcup\limits_{j \in J}{H_j^{I}}, \varphi_I \right], J \subset I.$$
Such a cube corresponds to a simplex in the descending link $\mathscr{L}$ of $[\Sigma,\id]$, namely the simplex associated to the edges from $[\Sigma, \mathrm{id}]$ to the vertices
$$\left[ \Sigma_I \backslash H^I_i, \varphi_I \right], \ i \in I.$$
We denote this simplex by $S_I$. 
Consider now the subcomplex $\mathscr{S}$ of the descending link $\mathscr{L}$ consisting of the union of all the simplices $S_I$. 

\begin{claim}\label{claim:subcplx}
The vertices of $\mathscr{S}$ correspond to the vertices in $V:=\{ [\Sigma_i \setminus H_i^i, \varphi_i] \mid i \in \mathscr{A} \}.$
\end{claim}

\noindent
Fix a set $I \subset \mathscr{A}$ of pairwise disjoint collections and an element $i \in I$. Notice that 
$$\left[ \Sigma_i \backslash H^i_i, \varphi_i \right]= \left[ \Sigma_I \backslash H^I_i, \varphi_I \right].$$
Indeed, $\varphi_i$ maps $\Sigma_i$ to $\Sigma$, $H^i_i$ to $D_i$, and is rigid outside $\Sigma$; and $\varphi_I$ maps $\Sigma_I$ to $\Sigma$, $H^I_i$ to $D_i$, is rigid outside $\Sigma$. So $\varphi_I^{-1} \varphi_i$ maps $\Sigma_i$ to $\Sigma_I$, $H^i_i$ to $H^I_i$, and is rigid outside $\Sigma_i$, Equivalently, $\varphi_I\varphi_i^{-1}$ sends $\Sigma_i \backslash H^i_i$ to $\Sigma_I \backslash H^I_i$ and is rigid outside $\Sigma_i \backslash H^i_i$. The desired equality follows. Thus, the vertices of $\mathscr{S}$ belong to $V$. Conversely, a vertex of $V$ is a simplex $S_J$ where $J$ is a singleton, concluding the proof of our claim.

\medskip \noindent
Our goal is now to show that $\mathscr{S}$ is a strict fundamental domain under the action of $\mathrm{Mod}(\Sigma)\leq \mathrm{stab}([\Sigma,\id])$ on the descending link $\mathscr{L}$. This is the content of the following two claims. 

\begin{claim}\label{claim:NoTwoSameOrbit}
No two simplices in $\mathscr{S}$ belong to the same $\mathrm{Mod}(\Sigma)$-orbit.
\end{claim}

\noindent
Two simplices in $\mathscr{S}$ can be written as $S_I$ and $S_J$ for some $I,J \subset \mathscr{A}$. Assume that there exists a $g\in\mathrm{Mod}(\Sigma)$ such that $g S_I=S_J$. If $i \in I$ is a frontier-arc, then $g$ must send the edge from $[\Sigma,\id]$ to $[\Sigma_I \backslash H^I_i, \varphi_I]$ (which corresponds to a vertex in $S_I$) to an edge from $[\Sigma,\id]$ to $[\Sigma_J \backslash H^J_j, \varphi_J]$ for some $j \in J$ (which is a vertex in $S_J$). We already know that $\varphi_J^{-1}g \varphi_I$ sends $\Sigma_I$ to $\Sigma_J$, and we deduce from the equality
$$\left[ \Sigma_I \backslash H^I_i, g \varphi_I \right] = \left[ \Sigma_J \backslash H^J_j, \varphi_J \right]$$
that it also sends $H^I_i \cap \mathrm{Fr}(\Sigma_I)$ to $H^J_j \cap \mathrm{Fr}(\Sigma_J)$. Because $\varphi^I_i$ sends $H^I_i \cap \mathrm{Fr}(\Sigma_I)$ to $i$ and since $g$ fixes pointwise $\partial \Sigma$, it follows that $\varphi_J^{-1}$ sends $i$ to $H^J_j \cap \mathrm{Fr}(\Sigma_J)$. But $\varphi_J$ sends $H^J_j \cap \mathrm{Fr}(\Sigma_J)$ to $j$, so necessarily $i=j \in J$. Thus, we have proved that $I \subset J$. By symmetry, we also have $J \subset I$, hence $S_I=S_J$. 

\begin{claim}\label{claim:HaveAtranslate}
Every simplex in $\mathscr{L}$ has a $\mathrm{Mod}(\Sigma)$-translate in $\mathscr{S}$.
\end{claim}

\noindent
Let us fix a simplex $S$ in $\mathscr{L}$. It corresponds to a cube $Q$, which is generated by the vertices
$$\left\{ \left[\Xi \cup \bigcup\limits_{j \in J} K_j, \psi \right] \mid J \subset \{1, \ldots, s\} \right\},$$
where $\Xi$ is an admissible subsurface, $K_1, \ldots, K_s$ are $2n$-gons of the rigid structure adjacent to $\Xi$, $\psi$ is an asymptotically rigid homeomorphism, and $[\Xi \cup  K_1 \cup \cdots \cup K_s, \psi]=[ \Sigma, \mathrm{id}]$. 

\medskip \noindent
Let $J \subset \mathscr{A}$ denote the set of pairwise disjoint collections corresponding to the intersections between the frontier of $\Sigma$ and the $\psi(K_i)$'s (which lie in $\Sigma$). So, for every $1 \leq i \leq s$, there exists some $j(i) \in J$ such that $\psi(K_i)$ and $D_{j(i)}$ have the same intersection with the frontier of $\Sigma$. Fix a $\beta \in \mathrm{Mod}(\Sigma)$ such that $\beta(\psi(K_i))= D_{j(i)}$ for every $1 \leq i \leq s$. The existence of such a $\beta$ follows from the fact that the discs $D_{j(i)}$ as well as the images $\psi(K_i)$ are disjoint. To summarize, we have:
\begin{itemize}
	\item The asymptotically rigid homeomorphism $\psi$ maps $\Xi \cup K_1 \cup \cdots \cup K_s$ to $\Sigma$, each of the $K_i$ to $\psi(K_i)$, and $\psi$ is rigid outside $\Xi \cup K_1 \cup \cdots \cup K_s$.
	\item The asymptotically rigid homeomorphism $\beta$ sends $\Sigma$ to $\Sigma$, each of the images $\psi(K_i)$ to  the disc $D_{j(i)}$, and it is rigid outside $\Sigma$.
	\item The asymptotically rigid homeomorphism $\varphi_I^{-1}$ maps $\Sigma$ to $\Sigma_I$, each of the discs $D_{j(i)}$ to the $2n$-gon $H^I_{j(i)}$, and it is rigid outside $\Sigma$.
\end{itemize}
This implies that $\varphi_I^{-1} \beta \psi$ maps $\Xi \cup K_1 \cup \cdots \cup K_s$ to $\Sigma$, each of the $K_i$ to $H^I_{j(i)}$, and it is rigid outside $\Xi \cup K_1 \cup \cdots \cup K_s$. Consequently,
$$\beta \cdot \left[ \Xi \cup K_1 \cup \cdots \cup K_{i-1} \cup K_{i+1} \cup \cdots \cup K_s, \psi \right] = \left[ \Sigma_I \backslash H^I_{j(i)}, \varphi_I \right]$$
for every $1 \leq i \leq s$. In other words, $\beta$ sends $S$ to $S_I$. This concludeds the proof of our claim.

\medskip \noindent
Thus, we have proved that $\mathscr{S}$ is a strict fundamental domain for the action of $\mathrm{Mod}(\Sigma)$ on $\mathscr{L}$. Next, notice that:

\begin{claim}\label{claim:StabModSimplex}
Fix a set $I \subset \mathscr{A}$ of pairwise disjoint collections. The stabiliser of $S_I$ in $\mathrm{Mod}(\Sigma)$ coincides with the subgroup $\mathrm{Mod}\left(\Sigma \backslash \bigcup\limits_{i \in I} D_i \right)$. 
\end{claim}

\noindent
Notice that, as a consequence of Claim \ref{claim:NoTwoSameOrbit}, the stabiliser of $S_I$ coincides with its pointwise stabiliser. If $g \in \mathrm{Mod}(\Sigma)$ stabilises $S_I$, then, for every $i \in I$, we have
$$\left[ \Sigma_I \backslash H^I_i, g \varphi_I \right] = g \cdot \left[ \Sigma_I \backslash H^I_i, \varphi_I \right]= \left[ \Sigma_I \backslash H^I_i, \varphi_I \right],$$
so $\varphi_I^{-1}g \varphi_I$ is isotopic to a homeomorphism that stabilises $\Sigma_I \backslash H^I_i$ and that is rigid outside. But we know that $\varphi_I^{-1}g \varphi_I$ first sends through $\varphi_I$ the subsurface $\Sigma_I$ to $\Sigma$, $H^I_i$ to $D_i$ and is rigid outside $\Sigma_I$; next it mixes $\Sigma$ by $g$, and is rigid outside $\Sigma$; finally it sends back $\Sigma$ to $\Sigma_I$ and $D_i$ to $H^I_i$ through $\varphi_I^{-1}$, and is rigid outside $\Sigma$. Therefore, $g$ has to preserve (up to isotopy) the topological disc $D_i$. This concludes the proof of our claim.
\begin{figure}
\begin{center}
\includegraphics[scale=0.3]{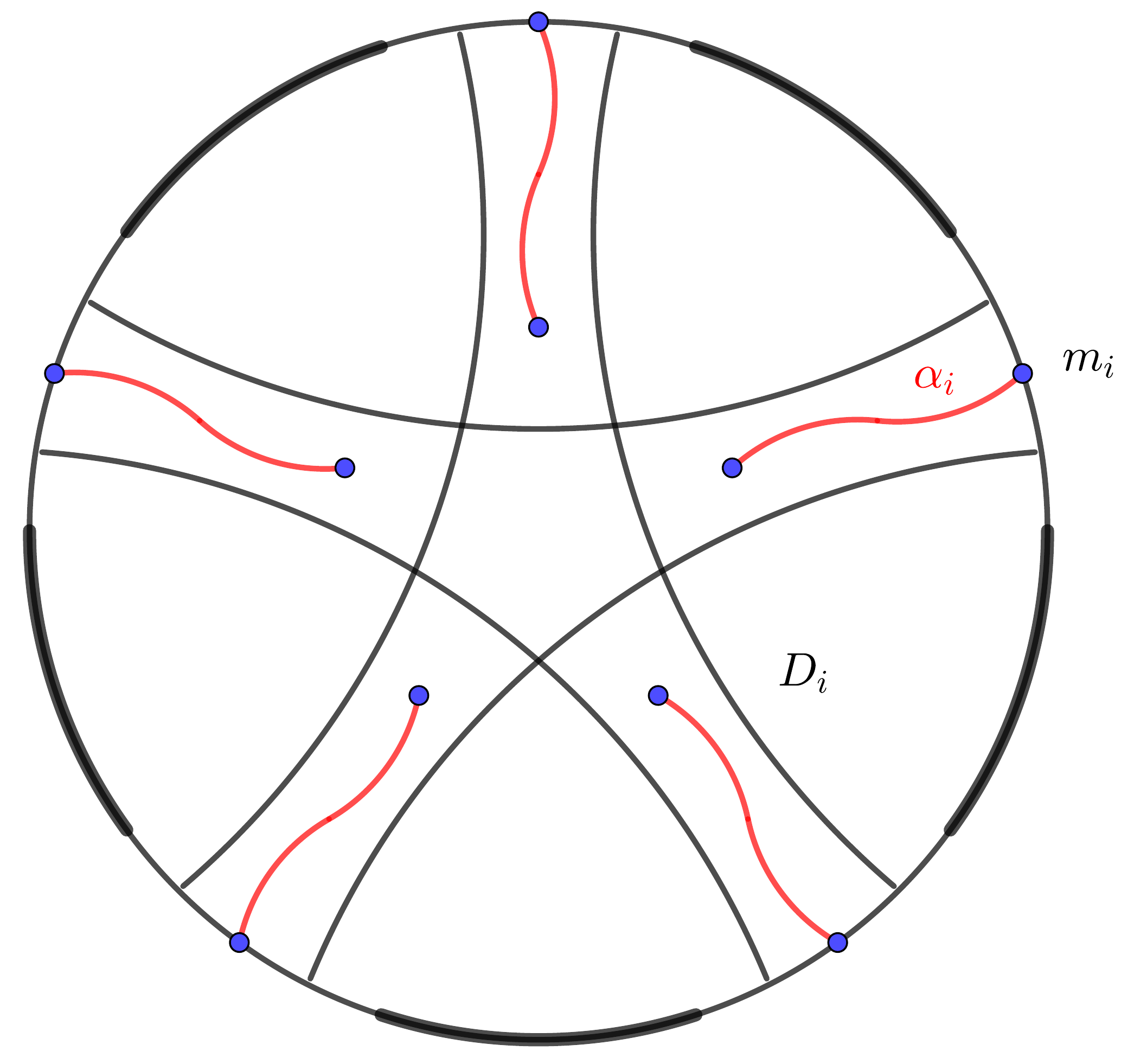}
\caption{Associating a marked point $m_i$ and an arc $\alpha_i$ to a disc $D_i$.}
\label{ArcLink}
\end{center}
\end{figure}

\medskip \noindent
For every $i \in \mathscr{A}$, fix a marked point $m_i \in \partial \Sigma$ that lies in the mid frontier-arc of $i$ if $n$ is odd or that lies between the two mid frontier-arcs of $i$ if $n$ is even. We identify the set of marked points $\{m_i \mid i \in \mathscr{A}\}$ with $\mathbb{Z}_c$ (where $c:=m+(k-1)(n-1)$ denotes the cardinality of $\mathscr{A}$) in such a way that the cyclic order induced by $\partial \Sigma$ coincides with the cyclic order induced by $\mathbb{Z}_c$. For every $i \in \mathscr{A}$, we also fix an arc $\alpha_i$ from the marked point $m_i$ to the unique puncture of the disc $D_i$. See Figure \ref{ArcLink}. Notice that $I \subset \mathscr{A}$ is a set of pairwise disjoint collections if and only if $\{D_i, i \in I\}$ is a collection of pairwise disjoint discs if and only if the arcs in $\{\alpha_i, i \in I\}$ are pairwise $(n-1)$-separated. By identifying $\Sigma$ with the disc $\mathbb{D}$ from the definition of $\mathfrak{C}(k,c,n-1)$, we deduce from the combination of Claims~\ref{claim:NoTwoSameOrbit},~\ref{claim:HaveAtranslate},~\ref{claim:StabModSimplex} with Lemma~\ref{lem:LinkTopo} that the map
$$g \cdot \left( [\Sigma,\mathrm{id}], \left[ \Sigma_{\{i\}} \backslash H_i^{\{i\}}, \varphi_{\{i\}} \right] \right) \mapsto g \cdot \alpha_i$$
induces an isomorphism $\mathscr{L} \to \mathfrak{C}(k,c,n-1)$.
\end{proof}

\begin{remark}
Notice that, as a consequence of \cite[Proposition II.12.20(1)]{BH} and of the proof of Proposition \ref{prop:deslink}, the descending links in $\mathscr{SC}(A_{n,m})$ can also be described as developments of simple complexes of braid groups. Interestingly, the underlying complexes turn out to coincide with the descending links in a CAT(0) cube complex on which the unbraided group $\amod_f(A_{n,m})$ acts. More precisely, if $n=1$, then the underlying complexes are $(m-1)$-simplices, which are isomorphic to the descending links in $\mathbb{R}^m$ on which $\amod_f(A_{n,m}) \simeq \mathbb{Z}^m \rtimes \mathbb{Z}_m$ acts; and if $n \geq 2$, then $\amod_f(A_{n,m}) \simeq T_{n,m}$ and the underlying complexes coincide with the descending links in the CAT(0) cube complex we can construct following \cite{MR2146639} by thinking of the Higmann-Thompson group $T_{n,m}$ as an \emph{annular diagram group}. 
\end{remark}

\subsection{Braided Higman-Thompson groups}\label{section:HigmanThompson}

\noindent
This section is dedicated to the proof of one of the two main results of this article, namely:


\ThompsonThm

\noindent
The strategy is to apply \emph{Morse theory} to the action of $\mathrm{br}(T_{n,m})$ on the spine $\mathscr{SC}(A_{n,m})$. The following criterion is a standard combination of Bestvina-Brady Morse theory \cite{Morse} and Brown's Criterion \cite{Brown}.

\begin{prop}\label{prop:Morse}
Let $G$ be a group acting cellularly on a contractible affine cell complex $X$ and let $f : X \to \mathbb{R}$ be a $G$-invariant Morse function. Assume that each sublevel set is $G$-cocompact, that the cell-stabilisers are of type $F_n$, and that, for every $k \geq 0$, there exists some $m \geq 0$ such that the descending link of every vertex $x \in X$ satisfying $f(x) \geq m$ is $k$-connected. Then $G$ is of type $F_n$. 
\end{prop}

\noindent
Consequently, in order to deduce Theorem \ref{PtolemyConnected} from Proposition \ref{prop:Morse}, we need to show that the complexes introduced in the previous section are highly connected. More precisely, our goal is to prove the following statement:

\begin{prop}\label{prop:ArcComplex}
Let $p\geq 2$, $q \geq 1$ and $r \geq 0$ be three integers. The complex $\mathfrak{C}(p,q,r)$ is $\left( \left\lfloor \frac{1}{3} \left( p + \left\lfloor \frac{q}{r+1} \right\rfloor \right) \right\rfloor -2 \right)$-connected. 
\end{prop}

\noindent
In fact, we are going to prove the proposition in a much more general framework. Indeed, in order to argue by induction, we need to introduce a more general family of complexes so that links of simplices still belong to this larger family. For the reader's convenience, let us fix the definitions of \emph{links} and \emph{stars} in simplicial complexes as used in this section and the next one.

\begin{definition}
Let $X$ be a simplicial complex and $\Delta \subset X$ a simplex. The \emph{link} of $\Delta$ is the subcomplex of $X$ that is the union of all the simplices that are disjoint from $\Delta$ but that span simplices with $\Delta$. The \emph{star} of $\Delta$, denoted by $\mathrm{star}(\Delta)$, is the union of all the simplices having $\Delta$ as a face.
\end{definition}

\noindent
The general framework in which we are going to prove Proposition \ref{prop:ArcComplex} is the following. Let $S$ be a punctured surface with boundary. Fix a set of punctures $P$, a set of marked points $M \subseteq \partial S$ and a symmetric relation $\sim$ on $M$. Here, we are interested in the simplicial complex $\mathfrak{R}=\mathfrak{R}(S,P,M,\sim)$ defined as follows: the vertices of $\mathfrak{R}$ are the isotopy classes of arcs connecting a point in $M$ to a point in $P$, and its simplices are collections of arcs that are pairwise disjoint and that start from marked points that are pairwise $\sim$-related. 

\medskip \noindent
Notice that, if $S$ is a disc with $p$ punctures, if $P$ is the set of all the punctures of $S$, if $M$ has cardinality $q$ and if $\sim$ is the $r$-separation, then $\mathfrak{R}$ coincides with the complex $\mathfrak{C}(p,q,r)$ from the previous section. Our main result about the connectedness properties of $\mathfrak{R}$ is the following:

\begin{prop}\label{prop:SimConnected}
Assume that $\#M \geq 1$ and $\#P \geq 2$. The complex $\mathfrak{R}(S,P,M,\sim)$ is $\left( \left\lfloor \frac{\#P + \min(\sim)}{3} \right\rfloor -2 \right)$-connected, where $\min(\sim)$ denotes the minimum size of a $\subseteq$-maximal collection of pairwise $\sim$-related marked points in $M$.
\end{prop}

\noindent
We record the following statement for future use during the proof of Proposition \ref{prop:SimConnected}.

\begin{lemma}\label{LEMMA39}
Let $Y$ be a compact $m$-dimensional combinatorial manifold. Let $X$ be a simplicial complex and assume that the link of every $k$-simplex in $X$ is $(m-k-2)$-connected. Let $\psi : Y \to X$ be a simplicial map whose restriction to $\partial Y$ is injective on simplices. Then, after possibly subdividing the simplicial structure of $Y$, $\psi$ is homotopic relative to $\partial Y$ to a map that is injective on simplices.
\end{lemma}

\noindent
This lemma appears in \cite[Lemma 3.9]{MR3545879} with the constant $m-2k-2$ instead of $m-k-2$. However, as it was pointed out to us by K.-U. Bux, there is a mistake in the proof of \cite[Lemma 3.9]{MR3545879}. (More precisely, at the very end of the first paragraph of the proof, the induction hypothesis does not apply to $\varphi : B \to \mathrm{link}(\psi(\sigma))$ because the assumptions only show that the link of a $d$-simplex in $\mathrm{link}(\psi(\sigma))$ is $(m-k-2d-3)$-connected instead of $(m-k-2d-2)$-connected.) Replacing the constant $m-2k-2$ with $m-k-2$ solves the problem, and Lemma \ref{LEMMA39} can be proved by reproducing the proof of \cite[Lemma 3.9]{MR3545879} word for word.

\begin{proof}[Proof of Proposition \ref{prop:SimConnected}.]
Our argument follows closely the proof of \cite[Theorem 3.10]{MR3545879}. We argue by induction over $\#M+\#P$. If $\#M+\#P=3$, then the statement is clear. From now on, we assume that $\#M+\#P > 3$. Fix a puncture $p \in P$ and a marked point $m \in M$, and let $\mathfrak{R}_0$ denote the subcomplex generated by the vertices corresponding to the arcs connecting a point in $M \backslash \{m\}$ to a point in $P \backslash \{p\}$. The first step of our argument is to show that we can work with $\mathfrak{R}_0$ instead of $\mathfrak{R}$.

\begin{claim}\label{claim:LinkInRzero1}
For every $k \geq 0$, the link of a $k$-simplex in $\mathfrak{R}_0$ is $\left( \left\lfloor \frac{\#P+\min(\sim)}{3} \right\rfloor - k-3 \right)$-connected. 
\end{claim}

\noindent
Let $x_0,\ldots, x_k$ denote the vertices of a $k$-simplex $\Delta$ in $\mathfrak{R}_0$. For every $0 \leq i \leq k$, $x_i$ is represented by an arc $\alpha_i$ connecting a point $n_i \in M \backslash \{m\}$ to a point $q_i \in P \backslash \{p\}$. Notice that the marked points $n_0, \ldots, n_k$ and the punctures $q_0, \ldots,q_k$ are pairwise distinct. By definition, the simplices in $\mathrm{link}(\Delta)$ correspond to the simplices in $\mathfrak{R}_0$ whose vertices are represented by arcs that are pairwse disjoint up to isotopy, that are disjoint from $\alpha_0,\ldots, \alpha_k$ up to isotopy, that start from pairwise $\sim$-related marked points, and that start from marked points $\sim$-related to $n_0, \ldots, n_k$. Consequently, the link of $\Delta$ is isomorphic to 
$$\mathfrak{R}(S\cup \{q_0, \ldots, q_k\}, P\backslash \{p, q_0, \ldots, q_k\}, M', \approx)$$
where $M'$ denotes the set of the elements in $M\backslash \{m,n_0, \ldots, n_k\}$ that are $\sim$-related to $n_0, \ldots, n_k$ and where $\approx$ denotes the restriction of $\sim$ to $M'$. By our induction hypothesis, we know that our link is $\left( \left\lfloor \frac{\#P-k-2+ \min(\approx)}{3} \right\rfloor -2 \right)$-connected. Now, fix a $\subseteq$-maximal collection in $M'$ of pairwise $\approx$-related points $\{y_1, \ldots, y_r\}$. Then either $\{y_1, \ldots, y_r, n_0, \ldots, n_k\}$ or $\{y_1, \ldots, y_r,m,n_0, \ldots, n_k\}$ is a $\subseteq$-maximal collection in $M$ of pairwise $\sim$-related points. Consequently, we have $\min(\sim) \leq r+k+2 = \min(\approx)+k+2$, which leads to the desired conclusion. 

\medskip \noindent
In the sequel, the previous argument will be used several times in different contexts.

\begin{claim}\label{claim:Rzero1}
The pair $(\mathfrak{R},\mathfrak{R}_0)$ is $\left( \left\lfloor \frac{\#P + \min(\sim)}{3} \right\rfloor -2 \right)$-connected, i.e. the inclusion $\mathfrak{R}_0 \hookrightarrow \mathfrak{R}$ induces an isomorphism on $\pi_i$ for $i< \left\lfloor \frac{\#P + \min(\sim)}{3} \right\rfloor -2$ and an epimorphism on $\pi_{\left\lfloor \frac{\#P + \min(\sim)}{3} \right\rfloor -2}$. 
\end{claim}

\noindent
Let $\mathfrak{R}_1$ denote the subcomplex generated by $\mathfrak{R}_0$ and the vertices, said of type $1$, corresponding to the arcs connecting $m$ to $p$. Because no two vertices of type $1$ are adjacent, $\mathfrak{R}_1$ is obtained from $\mathfrak{R}_0$ by gluing cones over the links of the vertices of type $1$. Arguing as in Claim \ref{claim:LinkInRzero1}, we show that such links are isomorphic to $\mathfrak{R}(S\cup \{p\},P\backslash \{p\}, M', \approx)$, where $M'$ denotes the set of the elements in $M \backslash \{m\}$ that are $\sim$-related to $m$ and where $\approx$ denotes the restriction of $\sim$ to $M'$. We know by induction that such a link is $\left( \left\lfloor \frac{\#P-1+\min(\approx)}{3} \right\rfloor -2 \right)$-connected and we show that the inequality $\min(\approx) \geq \min(\sim)-1$ holds, which allows us to deduce that the link under consideration is $\left( \left\lfloor \frac{\#P + \min(\sim)}{3} \right\rfloor -3 \right)$-connected. 

\medskip \noindent
Next, let $\mathfrak{R}_2$ denote the subcomplex generated by $\mathfrak{R}_1$ and the vertices, said of type $2$, corresponding to the arcs connecting $m$ to a point in $P \backslash \{p\}$. Because no two vertices of type $2$ are adjacent, $\mathfrak{R}_2$ is obtained from $\mathfrak{R}_1$ by gluing cones over the links of the vertices of type $2$. Arguing as in Claim \ref{claim:LinkInRzero1}, we show that such links are isomorphic to $\mathfrak{R}(S \cup \{q\}, P\backslash \{p,q\}, M', \approx)$, where $q$ is a puncture distinct from $p$, where $M'$ is the set of the elements in $M\backslash \{m\}$ that are $\sim$-related to $m$ and where $\approx$ denotes the restriction of $\sim$ to $M'$, and we show that they are $\left( \left\lfloor \frac{\#P + \min(\sim)}{3} \right\rfloor -3 \right)$-connected.

\medskip \noindent
Finally, let $\mathfrak{R}_3$ denote the subcomplex generated by $\mathfrak{R}_2$ and the vertices, said of type $3$, corresponding to the arcs connecting a point in $M \backslash \{m\}$ to $p$. Because no two vertices of type $3$ are adjacent, $\mathfrak{R}_3$ is obtained from $\mathfrak{R}_2$ by gluing cones over the links of the vertices of type $3$. Arguing as in Claim \ref{claim:LinkInRzero1}, we show that such links are isomorphic to $\mathfrak{R}(S\cup \{p\}, P\backslash \{p\}, M', \approx)$, where $n$ is a marked point distinct from $m$, where $M'$ is the set of the elements in $M\backslash \{n\}$ that are $\sim$-related to $n$ and where $\approx$ denotes the restriction of $\sim$ to $M'$, and that they are $\left( \left\lfloor \frac{\#P + \min(\sim)}{3} \right\rfloor -3 \right)$-connected. 

\medskip \noindent
Notice that a vertex in $\mathfrak{R}$ either belongs to $\mathfrak{R}_0$ or is of type $1$, $2$ or $3$, i.e. $\mathfrak{R}_3$ coincides with the entire complex $\mathfrak{R}$. Consequently, it follows from the previous paragraphs that $\mathfrak{R}$ can be obtained from $\mathfrak{R}_0$ by gluing cones over $\left( \left\lfloor \frac{\#P + \min(\sim)}{3} \right\rfloor -3 \right)$-connected subcomplexes, concluding the proof of Claim \ref{claim:Rzero1}.

\medskip \noindent
As a consequence of Claim \ref{claim:Rzero1}, it suffices to show that a map $\psi : \mathbb{S}^r \to \mathfrak{R}_0$ from a combinatorial sphere of dimension $r \leq \left\lfloor \frac{\#P + \min(\sim)}{3} \right\rfloor -2$ is homotopically trivial in $\mathfrak{R}$ in order to deduce that $\mathfrak{R}$ is $\left( \left\lfloor \frac{\#P + \min(\sim)}{3} \right\rfloor -2 \right)$-connected. By simplicial approximation, we may suppose without loss of generality that $\psi$ is simplicial. Also, as a consequence of Lemma \ref{LEMMA39}, which applies according to Claim \ref{claim:LinkInRzero1}, we may suppose without loss of generality that $\psi$ is injective on each simplex. 

\medskip \noindent
Fix an arc $\gamma$ from $m$ to $p$. We want to prove that $\psi$ can be homotoped so that its image lies in the star of $\gamma$. Since the latter is contractible, this will show that $\psi$ is homotopy trivial, as desired. 

\medskip \noindent
The arcs representing the vertices in the image of $\psi$ have their endpoints distinct from $p$ and $m$, but they may intersect $\gamma$. If there is no such intersection, then the vertices of the image of $\psi$ already lies in the star of $\gamma$. Consequently, the image of $\psi$ lies in the subcomplex generated by the star of $\gamma$, which coincides with the star of $\gamma$ itself because the link of $\gamma$ is \emph{flag} (i.e. every collection of pairwise adjacent vertices spans a simplex), so there is nothing to prove in this case. Otherwise, let $x \in \mathbb{S}^r$ be the vertex whose image is represented by the arc $\alpha$ that intersects $\gamma$ the closest to $p$. Fix a small disc $D \subseteq S$ containing $p$ such that $D \cap \alpha$ is a subarc contained in $\partial D$ and such that $D$ is disjoint from all the arcs representing the images under $\psi$ of the vertices of $\mathbb{S}^r$ distinct from $x$. Now let $\alpha'$ denote the arc obtained from $\alpha$ by replacing the subarc $\alpha \cap \partial D$ with $\partial D \backslash \alpha$. See Figure \ref{pushing}. Because $\psi$ is injective on simplices, the link of $x$ is sent in the link of $\psi(x)$ (which is represented by $\alpha$); and, by construction, this image also lies in the link of the vertex represented by $\alpha'$. Therefore, we can define a new map $\psi' : \mathbb{S}^r \to \mathfrak{R}_0$ by sending $x$ to the vertex represented by $\alpha'$ and by sending each vertex $y$ distinct from $x$ to $\psi(y)$. 
\begin{figure}
\begin{center}
\includegraphics[scale=0.32]{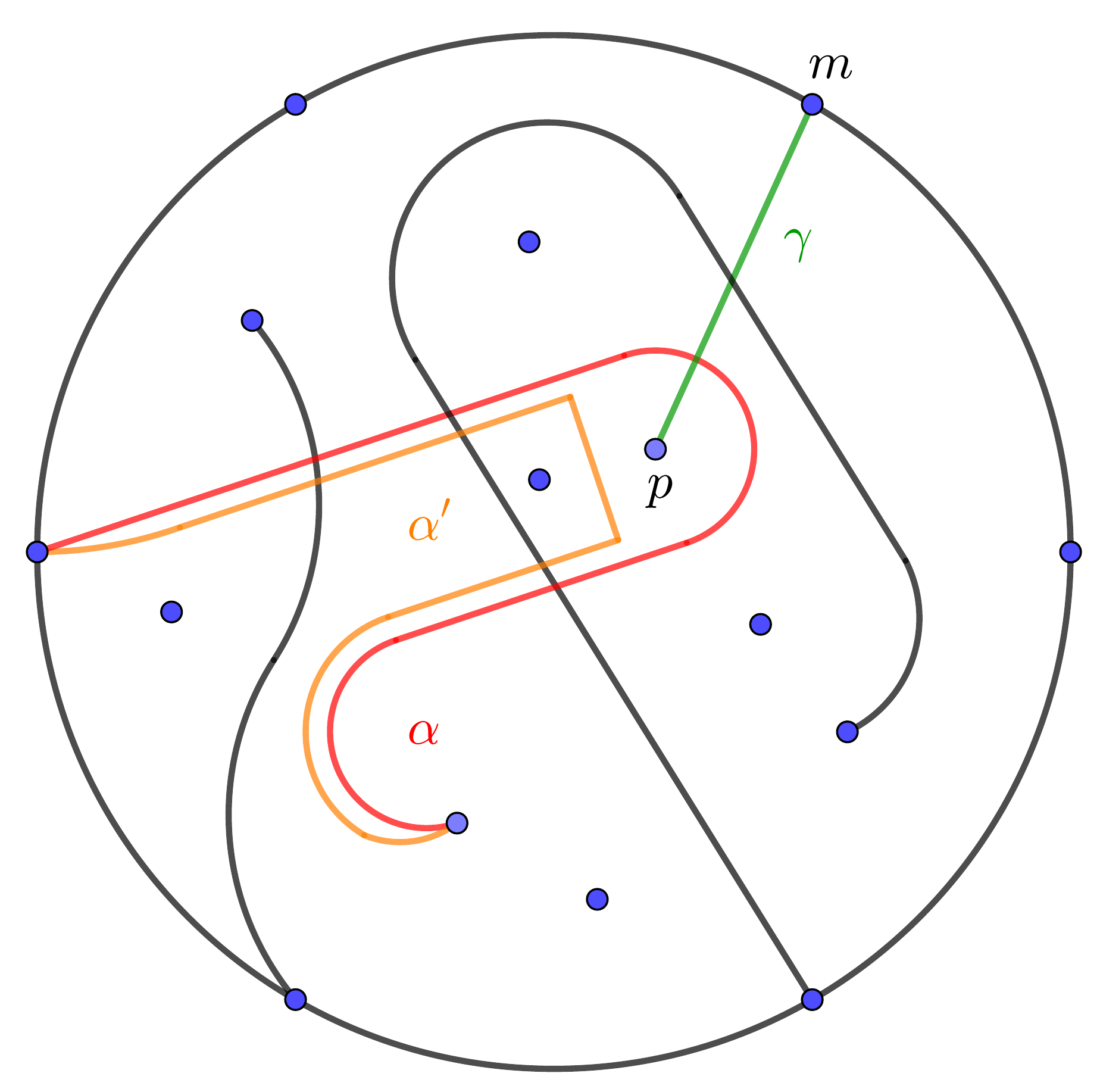}
\caption{Pushing $\alpha$ to $\alpha'$.}
\label{pushing}
\end{center}
\end{figure}

\medskip \noindent
We claim that $\psi$ and $\psi'$ are homotopy equivalent in $\mathfrak{R}$. Arguing as in Claim \ref{claim:LinkInRzero1}, we show that the intersection $L$ of the links in $\mathfrak{R}$ of the two vertices represented by $\alpha$ and $\alpha'$ is isomorphic to $\mathfrak{R}(S \cup \{p,q\}, P\backslash \{p,q\}, M', \approx)$, where $q \in P$ and $n \in M$ are the endpoints of $\alpha$, where $M'$ denotes the set of the elements in $M \backslash \{n\}$ that are $\sim$-related to $n$ and where $\approx$ denotes the restriction of $\sim$ to $M'$; and we show that $L$ is $\left( \left\lfloor \frac{\#P + \min(\sim)}{3} \right\rfloor -3 \right)$-connected. As a consequence, the common restriction $\mathbb{S}^{r-1} \to L$ of $\psi$ and $\psi'$ to $\mathrm{link}(x)$ is homotopically trivial, i.e. there exists a map $\varphi : \mathrm{star}(x) \to L$ such that $\varphi_{|\mathrm{link}(x)}$ coincides with the previous restriction. Because the star of $\psi(x)$ is contractible and that the image of $\varphi$ lies in the link of $\psi(x)$, we can homotope $\psi$ without modifying it outside the star of $x$ so that $\psi_{|\mathrm{star}(x)}= \varphi$. The same process applied to $\psi'$ leads to the same map, proving that $\psi$ and $\psi'$ are homotopy equivalent, as claimed.

\medskip \noindent
Notice that the total number of intersections between $\gamma$ and the arcs representing the images under $\psi'$ of the vertices in $\mathbb{S}^r$ is smaller than the total number of intersections between $\gamma$ and the arcs representing the images under $\psi$ of the vertices in $\mathbb{S}^r$. By iterating the argument, we construct a map $\mathbb{S}^r \to \mathfrak{R}_0$ that is homotopy equivalent to $\psi$ and whose image lies in the star of $\gamma$, as desired.

\medskip \noindent
Thus, we have proved that $\mathfrak{R}(S,P,M,\sim)$ is $\left( \left\lfloor \frac{\#P + \min(\sim)}{3} \right\rfloor -2 \right)$-connected, as desired.
\end{proof}

\begin{proof}[Proof of Proposition \ref{prop:ArcComplex}.]
Let $S$ be a disc with a set $P$ of $p$ punctures and a set $M$ of $q$ marked points in its boundary. Let $\sim$ denote the relation of being $r$-separated in $M$. Then $\mathfrak{R}(S,P,M,\sim)$ coincides with $\mathfrak{C}(p,q,r)$. By noticing that $\min(\sim)= \left\lfloor \frac{q}{r+1} \right\rfloor$, the desired conclusion follows from Proposition \ref{prop:SimConnected}.
\end{proof}

\begin{proof}[Proof of Theorem \ref{PtolemyConnected}.]
We want to apply Proposition \ref{prop:Morse} to the spine $\mathscr{SC}(A_{n,m})$ endowed with its height function. Notice that $\mathscr{SC}(A_{n,m})$ is contractible according to Proposition \ref{prop:HomotopyEqui} and Theorem \ref{thm:contractible}; that cell-stabilisers are finite extensions of braid groups according to Lemma \ref{lem:VertexStab}, so they are of type $F_\infty$; that, as a consequence of Proposition~\ref{prop:deslink} and Proposition \ref{prop:ArcComplex}, the descending link of a vertex is arbitrarily connected if its height is sufficiently high; and finally that Claim \ref{claim:CocompactLevel} below shows that sublevels are $\mathrm{br}T_{n,m}$-cocompact. Therefore, Proposition \ref{prop:Morse} proves that $\mathrm{br}T_{n,m}$ is of type $F_\infty$, as desired.

\begin{claim}\label{claim:CocompactLevel}
Fix two integers $n,m \geq 1$. For every $i \geq 1$, $\amod(A_{n,m})$ acts cocompactly on the subcomplex $X_i$ of $\mathscr{SC}(A_{n,m})$ generated by the vertices of height $\leq i$. 
\end{claim}

\noindent
Fix an $i \geq 1$ and let $C$ be a cube in $X_i$. So there exist an admissible subsurface $\Sigma$ containing the central polygon, a homeomorphism $\varphi$ and polygons $H_1, \ldots, H_k$ such that the vertices of $C$ are
$$\left\{ \left[ \Sigma \cup \bigcup\limits_{j \in I} H_j, \varphi \right] \mid I \subseteq \{1, \ldots, k\} \right\}.$$
Consequently, the vertices of the $\amod(A_{n,m})$-translate $\varphi^{-1} \cdot C$ are
$$\left\{ \left[ \Sigma \cup \bigcup\limits_{j \in I} H_j, \id \right] \mid I \subseteq \{1, \ldots, k\} \right\}.$$
But $\Sigma$ must be the union of at most $i$ polygons and it contains the central polygons, so there are only finitely many possibilities for $\Sigma$ and $H_1, \ldots, H_k$. Thus, $X_i$ contains only finitely many $\amod(A_{n,m})$-translates of cubes, concluding the proof of our claim.
\end{proof}

\begin{remark}
We do not expect the constant given by Proposition \ref{prop:ArcComplex} to be optimal. In fact, we think that $\mathfrak{R}(S,P,M,\sim)$ is homotopy equivalent to a bouquet of infinitely many $\left( \min \left( \#P, \min(\sim) \right)-1 \right)$-spheres. Proposition \ref{prop:BouquetSpheres} below proves this assertion in a particular case. However, the weaker information provided by Proposition \ref{prop:ArcComplex} is sufficient for our purpose, so we do not pursue this direction further.
\end{remark}

\subsection{Braided Houghton groups}\label{section:Houghton}

\noindent
This section is dedicated to the second main result of this article, namely:

\HoughtonThm

\noindent
Similarly to the unbraided Houghton groups, the strategy is to apply the following criterion:

\begin{prop}\label{prop:BrownHoughton}\emph{\cite[Corollary 3.3]{Brown}}
Let $G$ be a group acting on a contractible simplicial complex $X$ with cell-stabilisers that are finitely presented and of type $FP_\infty$. Fix a filtration $X_1 \subset X_2 \subset \cdots$ such that $G$ acts cocompactly on each $X_i$. If, up to homotopy, $X_{i+1}$ is obtained from $X_i$ by the adjunction of $n$-cells for every $i \geq 1$, then $G$ is of type $FP_{n-1}$ but not $FP_n$. Moreover, if $n \geq 3$, then $G$ is finitely presented.
\end{prop}

\noindent
More precisely, we are going to show that the descending links in the cube complex $\mathscr{SC}(A_{1,n})$, on which $\mathrm{br}H_n$, acts are homotopy equivalent to bouquets of $(n-1)$-spheres, which amounts to proving, as a consequence of Proposition \ref{prop:deslink}, that the complex $\mathfrak{C}(k,n,0)$ is $(n-2)$-connected (when $k$ is large enough) {by theorems of Hurewicz and Whitehead}. We already know from Proposition \ref{prop:ArcComplex} that this complex is $\left( \left\lfloor \frac{k+n}{3} \right\rfloor -2 \right)$-connected, but this is not sufficient. However, assuming that $r=0$, we will be able to reproduce the proof of Proposition \ref{prop:ArcComplex} almost word by word but with optimal constants. 

\medskip \noindent
Let $S$ be a punctured surface with boundary. Fix a set of punctures $P$ and a set of marked points $M \subset \partial S$. Here, we are interested in the simplicial complex $\mathfrak{R}(S,P,M)$ defined as follows: the vertices of $\mathfrak{R}(S,P,M)$ are the isotopy classes of arcs connecting a point in $M$ to a point in $P$, and its simplices are collections of arcs that are pairwise disjoint up to isotopy. Notice that, if $S$ is a disc with $p$ punctures, if $P$ is the set of all the punctures of $S$, and if $M$ has cardinality $q$, then $\mathfrak{R}(S,P,M)$ coincides with the complex $\mathfrak{C}(p,q,0)$.

\begin{prop}\label{prop:BouquetSpheres}
If $\#M \geq 1$ and $\#P \geq 2 \cdot \#M$, then $\mathfrak{R}(S,P,M)$ is homotopy equivalent to a bouquet of infinitely many $(\#M-1)$-spheres.
\end{prop}

\begin{proof}
First of all, we claim that $\mathfrak{R}(S,P,M)$ is $(\#M-2)$-connected. Our argument follows closely the proof of \cite[Theorem 3.10]{MR3545879}. We argue by induction over $\#M$. If $\#M=1$ the statement is clear and there is nothing to prove. So assume that $\#M \geq 2$. Fix a puncture $p \in P$ and a marked point $m \in M$, and let $\mathfrak{R}_0$ denote the subcomplex generated by the vertices corresponding to the arcs connecting a point in $M \backslash \{m\}$ to a point in $P \backslash \{p\}$.

\begin{claim}\label{claim:LinkInRzero}
For every $k \geq 0$, the link of a $k$-simplex in $\mathfrak{R}_0$ is $(\#M-k-4)$-connected.
\end{claim}

\noindent
Let $x_0,\ldots, x_k$ denote the vertices of a $k$-simplex $\Delta$ in $\mathfrak{R}_0$. For every $0 \leq i \leq k$, $x_i$ is represented by an arc $\alpha_i$ connecting a point $n_i \in M \backslash \{m\}$ to a point $q_i \in P \backslash \{p\}$. Notice that the marked points $n_0, \ldots, n_k$ and the punctures $q_0, \ldots,q_k$ are pairwise distinct.
By definition, the simplices in $\mathrm{link}(\Delta)$ correspond to the simplices in $\mathfrak{R}_0$ whose vertices are represented by arcs that are pairwse disjoint up to isotopy and that are disjoint from $\alpha_0,\ldots, \alpha_k$ up to isotopy.
Consequently, the link of $\Delta$ is isomorphic to
$$\mathfrak{R}(S\cup \{q_0, \ldots, q_k\}, P\backslash \{p, q_0, \ldots, q_k\}, M \backslash \{m,n_0, \ldots, n_k\}). $$
By induction, this complex is $(\#M-k-4)$-connected, as desired.

\begin{claim}\label{claim:Rzero}
The pair $(\mathfrak{R}(S,P,M),\mathfrak{R}_0)$ is $(\#M-2)$-connected, i.e. the inclusion $\mathfrak{R}_0 \hookrightarrow \mathfrak{R}(S,P,M)$ induces an isomorphism on $\pi_i$ for $i< \#M-2$ and an epimorphism on $\pi_{\#M-2}$. 
\end{claim}

\noindent
Let $\mathfrak{R}_1$ denote the subcomplex generated by $\mathfrak{R}_0$ and the vertices, said of type $1$, corresponding to the arcs connecting $m$ to $p$. Because no two vertices of type $1$ are adjacent, $\mathfrak{R}_1$ is obtained from $\mathfrak{R}_0$ by gluing cones over the links of the vertices of type $1$. Notice that such links are isomorphic to $\mathfrak{R}(S\cup \{p\},P\backslash \{p\}, M\backslash \{m\})$, and so are $(\#M-3)$-connected by induction. 

\medskip \noindent
Next, let $\mathfrak{R}_2$ denote the subcomplex generated by $\mathfrak{R}_1$ and the vertices, said of type $2$, corresponding to the arcs connecting $m$ to a point in $P \backslash \{p\}$. Because no two vertices of type $2$ are adjacent, $\mathfrak{R}_2$ is obtained from $\mathfrak{R}_1$ by gluing cones over the links of the vertices of type $2$. Notice that such links are isomorphic to $\mathfrak{R}(S \cup \{q\}, P\backslash \{p,q\}, M \backslash \{m\})$ where $q$ is a puncture distinct from $p$, and so they are $(\#M-3)$-connected by induction.

\medskip \noindent
Finally, let $\mathfrak{R}_3$ denote the subcomplex generated by $\mathfrak{R}_2$ and the vertices, said of type $3$, corresponding to the arcs connecting a point in $M \backslash \{m\}$ to $p$. Because no two vertices of type $3$ are adjacent, $\mathfrak{R}_3$ is obtained from $\mathfrak{R}_2$ by gluing cones over the links of the vertices of type $3$. Notice that such links are isomorphic to $\mathfrak{R}(S\cup \{p\}, P\backslash \{p\}, M\backslash \{n\})$ where $n$ is a marked point distinct from $m$, and so they are $(\#M-3)$-connected by induction. 

\medskip \noindent
Notice that a vertex in $\mathfrak{R}(S,P,M)$ either belongs to $\mathfrak{R}_0$ or it is of type $1$, $2$ or $3$, i.e. $\mathfrak{R}_3$ coincides with the entire complex $\mathfrak{R}(S,P,M)$. Consequently, it follows from the previous paragraphs that $\mathfrak{R}(S,P,M)$ can be obtained from $\mathfrak{R}_0$ by gluing cones over $(\#M-3)$-connected subcomplexes. This concludes the proof of Claim \ref{claim:Rzero}.

\medskip \noindent
As a consequence of Claim \ref{claim:Rzero}, it suffices to show that a map $\psi : \mathbb{S}^r \to \mathfrak{R}_0$ from a combinatorial sphere of dimension $r \leq \#M-2$ is homotopically trivial in $\mathfrak{R}(S,P,M)$ in order to deduce that $\mathfrak{R}(S,P,M)$ is $(\#M-2)$-connected. By simplicial approximation, we may suppose without loss of generality that $\psi$ is simplicial. Also, as a consequence of Lemma \ref{LEMMA39}, which applies according to Claim \ref{claim:LinkInRzero}, we may suppose without loss of generality that $\psi$ is injective on each simplex. 

\medskip \noindent
Fix an arc $\gamma$ from $m$ to $p$. We want to prove that $\psi$ can be homotoped so that its image lies in the star of $\gamma$  in $\mathfrak{R}(S,P,M)$. Since the star of a vertex is contractible, this will show that $\psi$ is homotopy trivial, as desired. 

\medskip \noindent
The arcs representing the vertices in the image of $\psi$ have their endpoints distinct from $p$ and $m$, but they may intersect $\gamma$. If there is no such intersection, then the vertices of the image of $\psi$ already lies in the star of $\gamma$. Consequently, the image of $\psi$ lies in the subcomplex generated by the star of $\gamma$, which coincides with the star of $\gamma$ itself because the link of $\gamma$ is \emph{flag} (i.e. every collection of pairwise adjacent vertices spans a simplex), so there is nothing to prove in this case. Otherwise, let $x \in \mathbb{S}^r$ be the vertex whose image is represented by the arc $\alpha$ that intersects $\gamma$ the closest to $p$. Fix a small disc $D \subset S$ containing $p$ such that $D \cap \alpha$ is a subarc contained in $\partial D$ and such that $D$ is disjoint from all the arcs representing the images under $\psi$ of the vertices of $\mathbb{S}^r$ distinct from $x$. Now let $\alpha'$ denote the arc obtained from $\alpha$ by replacing the subarc $\alpha \cap \partial D$ with $\partial D \backslash \alpha$. See Figure~\ref{pushing}. Because $\psi$ is injective on simplices, the link of $x$ is sent in the link of $\psi(x)$ (which is represented by $\alpha$); and, by construction, this image also lies in the link of the vertex represented by $\alpha'$. Therefore, we can define a new map $\psi' : \mathbb{S}^r \to \mathfrak{R}_0$ by sending $x$ to the vertex represented by $\alpha'$ and by sending each vertex $y$ distinct from $x$ to $\psi(y)$. 

\medskip \noindent
We claim that $\psi$ and $\psi'$ are homotopy equivalent in $\mathfrak{R}(S,P,M)$. Notice that the intersection $L$ of the links in $\mathfrak{R}(S,P,M)$ of the vertices represented by $\alpha$ and $\alpha'$ is isomorphic to $\mathfrak{R}(S \cup \{p,q\}, P\backslash \{p,q\}, M \backslash \{n\})$ where $q \in P$ and $n \in M$ are the endpoints of $\alpha$. By induction, the intersection is therefore $(\#M-3)$-connected. As a consequence, the common restriction $\mathbb{S}^{r-1} \to L$ of $\psi$ and $\psi'$ to $\mathrm{link}(x)$ is homotopically trivial, i.e. there exists a map $\varphi : \mathrm{star}(x) \to L$ such that $\varphi_{|\mathrm{link}(x)}$ coincides with the previous restriction. Because the star of $\psi(x)$ is contractible and the image of $\varphi$ lies in the link of $\psi(x)$, we can homotope $\psi$ without modifying it outside the star of $x$ so that $\psi_{|\mathrm{star}(x)}= \varphi$. The same process applied to $\psi'$ leads to the same map, proving that $\psi$ and $\psi'$ are homotopy equivalent, as claimed.

\medskip \noindent
Notice that the total number of intersections between $\gamma$ and the arcs representing the images under $\psi'$ of the vertices in $\mathbb{S}^r$ is smaller than the total number of intersections between $\gamma$ and the arcs representing the images under $\psi$ of the vertices in $\mathbb{S}^r$. By iterating the argument, we construct a map $\mathbb{S}^r \to \mathfrak{R}_0$ that is homotopy equivalent to $\psi$ and whose image lies in the star of $\gamma$, as desired.

\medskip \noindent
Thus, we have proved that $\mathfrak{R}(S,P,M)$ is $(\#M-2)$-connected. Notice that $\mathfrak{R}(S,P,M)$ has dimension $\#M-1$. But it is well-known that an $(n-1)$-connected CW complex of dimension $n$ is homotopy equivalent to a bouquet of $n$-spheres, so it follows that $\mathfrak{R}(S,P,M)$ is homotopy equivalent to a bouquet of $(\#M-1)$-spheres.

\medskip \noindent
It remains to show that this bouquet contains infinitely many spheres. Because the complex $\mathfrak{R}(S,P,M)$ is $(\#M-1)$-dimensional and $(\#M-2)$-connected, its $(\#M-1)$th homotopy group is isomorphic to its $(\#M-1)$th homology group, which is isomorphic to the free abelian group of $(\#M-1)$-chains. Therefore, in order to conclude, it suffices to construct infinitely many non-trivial $(\#M-1)$-chains using pairwise disjoint collections of simplices in $\mathfrak{R}(S,P,M)$.

\medskip \noindent
Write $M=\{m_1, \ldots, m_k\}$. Because $\#P \geq 2 \#M$, we can assign to each marked point $m_i$ two punctures $p_i,q_i$ such that $p_1,q_1, \ldots, p_k,q_k$ are pairwise distinct. For every $1 \leq i \leq k$, let $\alpha^i_1,\alpha^i_2,\ldots$ be a sequence of arcs connecting $m_i$ to $p_i$ and supported in a topological disc containing $m_i$, $p_i$ and $q_i$ among the marked points and punctures. We choose our arcs such that:
\begin{itemize}
	\item for every $1 \leq i \leq k$, the arcs $\alpha^i_1, \alpha^i_2, \ldots$ are pairwise non-isotopic;
	\item for all distinct $1 \leq i,j \leq k$ and for all $r,s \geq 1$, the arcs $\alpha^i_r$ and $\alpha^j_s$ are disjoint.
\end{itemize}
For all $j \geq 1$, let $S_j$ denote the subcomplex generated by the vertices represented by the arcs $\{\alpha^i_j, \alpha^i_{j+1} \mid 1 \leq i \leq k\}$. Notice that, for every $1 \leq t \leq k-1$, the subcomplex associated to $\{\alpha^i_j,\alpha^i_{j+1} \mid 1 \leq i \leq t+1\}$ coincides with the suspension of the subcomplex associated to $\{\alpha^i_j, \alpha^i_{j+1} \mid 1 \leq i \leq t\}$,  so $S_j$ is a triangulation of a $(k-1)$-sphere. Consequently, $S_j$ can be thought of as a $(k-1)$-chain, and it is non-trivial as a sum of pairwise distinct simplices. Moreover, for every $j \geq 1$, $S_j$ and $S_{j+2}$ are disjoint, hence infinitely many non-trivial and pairwise disjoint $(k-1)$-chains, as desired.
\end{proof}

\noindent
As an immediate consequence of Proposition \ref{prop:BouquetSpheres}, we get:

\begin{cor}\label{cor:ArcComplex}
For all $q \geq 1$ and $p \geq 2q$, the complex $\mathfrak{C}(p,q,0)$ is homotopy equivalent to a bouquet of infinitely many $(q-1)$-spheres.
\end{cor}

\begin{proof}[Proof of Theorem \ref{thm:brHfiniteness}.]
Fix an $n \geq 1$ and, for every $i \geq 1$, let $X_i$ denote the subcomplex of $\mathscr{SC}(R_n)$ generated by the vertices of height $\leq 2n+i$. Recall from Proposition \ref{prop:HomotopyEqui} and Theorem \ref{thm:contractible} that $\mathscr{SC}(A_{n,m})$ is contractible. Moreover, we know from Claim \ref{claim:CocompactLevel} that each $X_i$ is $\amod(A_{1,n})$-cocompact, and a fortiori $\mathrm{br}H_n$-cocompact since $\mathrm{br}H_n$ has finite index in $\amod(A_{1,n})$. Next, notice that, for every $i \geq 1$, $X_{i+1}$ is obtained from $X_i$ by gluing cones over the descending links of the vertices of height $i+1$ \cite[Lemma 2.5]{Morse}. According to Proposition \ref{prop:deslink} and Corollary \ref{cor:ArcComplex}, these links are homotopy equivalent to bouquets of $(n-1)$-spheres. As a consequence, up to homotopy, $X_{i+1}$ is obtained from $X_i$ by adjunctions of $n$-cells. Now, Theorem \ref{thm:brHfiniteness} follows from Proposition \ref{prop:BrownHoughton}. 
\end{proof}

\addcontentsline{toc}{section}{References}

\bibliographystyle{alpha}
{\footnotesize\bibliography{BraidedThompson}}

\Address

\end{document}